\documentclass[9pt,shortpaper,twoside,web]{IEEEtran}
\usepackage{cite}
\usepackage{amsmath,amssymb,amsfonts}
\usepackage{algorithmic}
\usepackage{graphicx}
\usepackage{textcomp}

\usepackage{epsfig} 
\usepackage{comment}
\usepackage{bm}
\usepackage{color}
\newtheorem{problem}{Problem}

\newtheorem{theorem}{Theorem}

\newtheorem{rmk}{Remark}
\newtheorem{assumption}{Assumption}
\newtheorem{conjecture}{Conjecture}

\newcommand{\argmax}{\mathop{\rm arg~max~}\limits}
\newcommand{\esssup}{\mathop{\rm ess~sup}\limits}
\newcommand{\sgn}{\mathrm{sgn}}

\newcommand{\red}[1]{{\color{black}{#1}}}
\newcommand{\blue}[1]{{\color{black}{#1}}}

\def\BibTeX{{\rm B\kern-.05em{\sc i\kern-.025em b}\kern-.08em
    T\kern-.1667em\lower.7ex\hbox{E}\kern-.125emX}}
\begin{document}
\title{Multiple sparsity constrained control node scheduling with application to rebalancing of mobility networks}
\author{Takuya Ikeda, \IEEEmembership{Member, IEEE}, 
	Kazunori Sakurama, \IEEEmembership{Member, IEEE}, and 
	Kenji Kashima, \IEEEmembership{Senior Member, IEEE}
\thanks{This work was partially supported by the joint project of Kyoto University and Toyota Motor
Corporation, titled ``Advanced Mathematical Science for Mobility Society''.}
\thanks{T. Ikeda is with 
		Faculty of Environmental Engineering, The University of Kitakyushu, Fukuoka, 808-0135, Japan
		(e-mail: t-ikeda@kitakyu-u.ac.jp).}
\thanks{K. Sakurama and K. Kashima are with
		Graduate School of Informatics, Kyoto University, Kyoto, 606-8501, Japan
		(e-mail: \{sakurama,kk\}@i.kyoto-u.ac.jp).}}
\maketitle

\begin{abstract}
This paper treats an optimal scheduling problem of control nodes in networked systems.
We newly introduce both the $L^0$ and $\ell^0$ constraints on control inputs to extract a time-varying small number of effective control nodes.
As the cost function, we adopt the trace of the controllability Gramian to reduce the required control energy.
Since the formulated optimization problem is combinatorial, we introduce a convex relaxation problem for its computational tractability.
After a reformulation of the problem into an optimal control problem to which Pontryagin's maximum principle is applicable,
we give a sufficient condition under which the relaxed problem gives a solution of the main problem.
Finally, the proposed method is  applied to a rebalancing problem of a mobility network.
\end{abstract}

\begin{IEEEkeywords}
convex optimization, networked systems, optimal control, sparse control
\end{IEEEkeywords}

\section{Introduction}\label{sec:introduction}
Nowadays, control system designs that incorporate a notion of sparsity have attracted a lot of attention in the control community. 
Such an approach finds essential information that gives a significant impact to the system of interest,
and it plays an important role in many occasions in large-scale networked systems such as {\em control node selection} tackled in this paper. 
There are mainly two types of penalty costs to enhance the sparsity.
The first one is the $\ell^0$ norm, which is defined as the number of non-zero components. 
This cost is widely used in sparse modeling motivated by the success of compressed sensing, 
and most of the related works in control systems also adopt this type.
The second one is the $L^0$ norm, which is defined as the length of the support. 
This is an extended version of the $\ell^0$ norm for functional spaces,
and it seems to appear in relatively recent works, e.g.,~\cite{ItoKun14,NagQueNes16,ITKKTAC18,KumSriCha19,IwaIzaKas19}.
However, it should be emphasized that 
optimization problems involving both of the $\ell^0$ norm and the $L^0$ norm
have not yet been investigated in the area of sparse optimization, to the best of our knowledge.

The purpose of the node selection problem is to identify the set of nodes that should receive exogenous control inputs 
so that the overall system of interest is effectively guided. 
The selected nodes are called {\it control nodes}.
This selection problem naturally arises in large-scale network systems due to physical or financial reasons.
In recent works, control nodes are chosen based on a metric of controllability. 
For example,
the work~\cite{Ols14} considers the minimum set of control nodes that ensures the classical controllability in~\cite{Kal63};
the work~\cite{AssKhaLiPre15} considers the structural controllability;
the works~\cite{PasZamBul14,SumCorLyg16} introduce quantities that evaluate how much the system is easy to control,
such as the trace of the controllability Gramian.

While the works above investigate the selection problem in which the set of control nodes is fixed over the time, 
more recent works alternatively consider time-varying control node selection, which is also referred to as {\em control node scheduling}. 
The node scheduling problem finds {\em not only which but also when nodes should become control nodes}, 
and hence it seems more challenging and efficient for achieving high control performance.
Indeed, the authors in~\cite{ZhaPasCor16,Noz17} consider the scheduling problem for discrete-time systems 
and show its effectiveness over the time-invariant control node selection.
Mathematically, all of the aforementioned works consider $\ell^0$ constrained optimization problems,
in which the number of nodes selected at the same time is constrained.
On the other hand, in networked systems
it is also important to effectively compress control signals 
and reduce communication traffic.
To achieve this, it is desirable to find the best time duration over which controllers should become active.
Then, we considered the $L^0$ constraint on control inputs
and formulated a node scheduling problem for continuous-time systems in~\cite{ITKKLCSS18}.
This scheduling problem is furthermore analyzed in~\cite{Ols19},
which provides an explicit formula of the optimal solutions 
and shows that the solutions are obtained by a greedy algorithm.
However, these two works on continuous-time systems mainly consider the $L^0$ control cost,
and the resulting number of activated control nodes at each time instance
(i.e., the $\ell^0$ control cost) is not taken into account.

In view of this, this paper newly considers an optimal node scheduling problem 
under the $L^0$ and $\ell^0$ constraints.
By introducing these two constraints, 
we can find a time-varying small number of control nodes while reducing the support of control inputs.
As the network controllability, we adopt the trace of the controllability Gramian.
This quantity is closely related to the average energy required to steer the system in all directions in the state space~\cite{MulWeb72}.
The formulated problem thus includes a combinatorial structure caused by the $L^0$ and $\ell^0$ norms.
To circumvent this, we introduce a convex relaxation problem 
and establish a condition for the main problem to be exactly solved via the convex optimization. 
For the analysis, we transform the convex relaxation problem to an optimal control problem to which Pontryagin's maximum principle is applicable.
Unlike the previous formulation in~\cite{ITKKLCSS18,Ols19},
our maximum condition in the principle does not boil down to component-wise calculation 
and the transversality condition is needed to show the equivalence 
between the main problem and the relaxation problem.

To demonstrate the practicability of the developed method, we address a rebalancing problem of mobility networks with one-way trips, e.g., car- and bike-sharing systems~\cite{Illgen2019,Ferrero2018}. In such a system, a customer picks up a vehicle in a station and can return it in another station, which enhances the usability of sharing systems but causes a problem of uneven distribution of vehicles. To maintain this system, rebalancing vehicles is required, which should be as infrequent as possible for the reduction of staff cost. This problem is formulated as  optimization problems to solve with conventional techniques of optimization and control engineering \cite{Zhang2019,Calafiore2019}.
Although in the existing research the amount of rebalanced vehicles has been taken into account, the frequency of rebalancing should be really considered to reduce staff cost.
In this paper, we show how to apply the developed method to attain infrequent rebalancing, namely sparse rebalancing, in the mobility network systems,
and illustrate its effectiveness through numerical examples.

The remainder of this paper is organized as follows: 
Section~\ref{sec:math} provides mathematical preliminaries.
Section~\ref{sec:formulation} formulates our node scheduling problem. 
Section~\ref{sec:analysis} introduces a convex relaxation problem and gives a sufficient condition for the main problem to boil down to the convex optimization.
A numerical example of the proposed node scheduling is also illustrated.
Section~\ref{sec:application_to_rebalancing} extends the proposed method to a rebalancing problem of mobility networks.
Section~\ref{sec:conclusion} offers concluding remarks.

\section{Mathematical Preliminaries}
\label{sec:math}

This section reviews notation that will be used throughout the paper.

We denote the set of all positive integers by $\mathbb{N}$
and the set of all real numbers by $\mathbb{R}$.
Let $m\in\mathbb{N}$ and $\Omega\subset\mathbb{R}$.
For a vector $a=[a_1, a_2, \dots, a_m]^{\top}\in\mathbb{R}^{m}$, 
$\mbox{diag}(a)$ denotes the diagonal matrix whose $(i, i)$-component is given by $a_i$,
and $a\in \Omega^{m}$ means $a_i\in \Omega$ for all $i$.
The {\em $\ell^0$ norm} and {\em $\ell^1$ norm} of $a$ are defined by
$\|a\|_{\ell^0} \triangleq \#\{i\in\{1, 2, \dots, m\}: a_i\neq0\}$ and
$\|a\|_{\ell^1} \triangleq \sum_{i=1}^{m} |a_i|,$
where $\#$ returns the number of elements of a set.
We denote the Euclidean norm by 
$\|a\|\triangleq (\sum_{i=1}^{m}a_{i}^{2})^{1/2}$.
Let $N\in\mathbb{N}$, $N_1\in\mathbb{N}$, and $N_2\in\mathbb{N}$.
For any $M\in\mathbb{R}^{N\times N}$, $\mathrm{Tr} M$ denotes the trace of $M$.
For any $M\in\mathbb{R}^{N_1\times N_2}$, $M^{\top}$ denotes the transpose of $M$.
Let $S$ be a closed subset of $\mathbb{R}^m$ and $a\in S$.
A vector $\xi\in\mathbb{R}^m$ is a proximal normal to the set $S$ at the point $a$
if and only if there exists a constant $\sigma\geq0$ such that
$\xi^\top(b-a) \leq \sigma \|b-a\|^2$ for all $b\in S$.
The {\em proximal normal cone} to $S$ at $a$ is defined as the set of all such $\xi$, 
which is denoted by $N_S^P(a)$.
We denote the {\em limiting normal cone} to $S$ at $a$ by $N_S^L(a)$, i.e.,
$N_S^L(a) \triangleq \{\xi = \lim_{i\to\infty} \xi_i: \xi_i \in N_S^P(a_i), a_i \to a, a_i \in S\}.$
Let $T>0$. 
We define the {\em $L^0$, $L^1$, $L^\infty$ norms} 
of a measurable function $v(t)=[v_1(t), v_2(t), \dots, v_m(t)]^{\top}\in{\mathbb{R}}^m$ on $[0, T]$ by
\begin{align*}
	&\|v\|_{L^0} \triangleq \sum_{j=1}^{m}\mu_{L}(\{t\in[0, T]: v_j(t)\neq0\}),\\
	&\|v\|_{L^1} \triangleq \sum_{j=1}^{m}\int_{0}^{T} |v_j(t)| dt,\\
	&\|v\|_{L^\infty} \triangleq \max_{1 \leq j \leq m} \esssup_{0\leq t\leq T} |v_j(t)|,	
\end{align*}
where $\mu_{L}$ is the Lebesgue measure on ${\mathbb{R}}$,
and $\esssup$ denotes the essential supremum defined by
\[
	\esssup_{0\leq t\leq T} |v_j(t)| \triangleq 
	\inf \left\{a\in\mathbb{R}: \mu_L\left(\left\{t\in[0, T]: |v_j(t)|>a\right\}\right)=0\right\}.
\]
We denote the set of all functions $v$ with $\| v \|_{p}<\infty$ by $L^p$, $p\in\{1,\infty\}$.
\red{We call a vector-valued function with absolutely continuous components {\em arc}~\cite[p.~255]{Cla13}.}

\section{Problem Formulation}
\label{sec:formulation}

\subsection{System Description}
\label{subsec:controllability}

Let us consider a network model consisting of $n$ nodes and define the overall system by
\begin{align}
\begin{split}
&\dot{x}(t)=Ax(t)+BV(t)u(t), \quad 0 \leq t \leq T,\\
&V(t)\triangleq \mathrm{diag}(v(t)),
\label{eq:system2}
\end{split}
\end{align}
where 
$x(t)=[x_1(t), x_2(t),\dots,x_n(t)]^\top\in{\mathbb{R}}^n$ is the state vector consisting of $n$ nodes, where $x_{i}(t)$ is the state of the $i$-th node at time $t$; 
$u(t)\in{\mathbb{R}}^m$ is the exogenous control input that influences the network dynamics;
$A\in{\mathbb{R}}^{n \times n}$ is the dynamics matrix that represents the information flow among nodes;
$B=[b_1,b_2,\dots,b_m]\in{\mathbb{R}}^{n\times m}$ is a constant matrix that represents candidates of control nodes;
$v(t)\in\{0, 1\}^{m}$ represents the activation schedule of the control input $u(t)$;
$T>0$ is the final time of control. 
In this setting, the control input $u_j(t)$, the $j$-th component of $u(t)$, is able to affect the system through the vector $b_j$ at time $t$ if and only if $v_j(t)=1$,
and the nodes that receive the inputs are called {\em control nodes}.
In other words, control node scheduling problem seeks an optimal variable $v(t)$ over $[0, T]$ based on a given cost function and some constraints.

\subsection{Main Problem}
\label{subsec:problem}
This paper is interested in the controllability performance as the cost function.
The performance is related to the quantity of the required control energy, 
for which a number of metrics have been proposed;
see e.g.~\cite{PasZamBul14,SumCorLyg16}.
Among them, this paper adopts the trace of the controllability Gramian,
which is a metric to approximate the network average controllability~\cite{MulWeb72}.
In addition, this paper introduces the $L^0$ and $\ell^0$ constraints on the control input.
The $L^0$ constraint limits the time duration where the control becomes active,
and the $\ell^0$ constraint limits the number of control nodes at each time.
Thus, the main problem of this paper is defined as follows:
\begin{problem}
\label{prob:L0}
Given $A\in\mathbb{R}^{n\times n}$, $B\in\mathbb{R}^{n\times m}$, $T>0$, $\beta\in\{1,2,\cdots,m-1\}$, and $\alpha_j\in(0, T]$, $j=1,2,\dots,m$,
find a time-varying matrix 
$V(t)\triangleq\mathrm{diag}(v(t))$,
$v(t)\triangleq[v_1(t), v_2(t), \dots, v_m(t)]^{\top}$,
that solves
\begin{equation*}
\begin{aligned}
	& \underset{v}{\text{maximize}}
	& & J(v)\triangleq\mathrm{Tr} \int_{0}^{T} e^{At}BV(t)V(t)^{\top} B^{\top}e^{A^{\top}t}dt\\
	& \text{subject to}
	& & v(t) \in \{0, 1\}^m \quad \forall t\in[0, T], \\
	& & & \|v_j\|_{L^0} \leq\alpha_j \quad \forall j \in\{1,2,\dots,m\},\\
	& & & \|v(t)\|_{\ell^0} \leq\beta \quad \forall t\in[0, T].
\end{aligned}
\end{equation*}
\end{problem}
In this paper, we will show that Problem~\ref{prob:L0} is exactly solved via a convex optimization problem.
Note that given two optimization problems are said to be {\em equivalent} 
if the set of all optimal solutions coincides. 

\begin{rmk}
Note that, from a property of the trace operator
we have 
\[
	J(v) = \mathrm{Tr} \int_{0}^{T} B^{\top}e^{A^{\top}t} e^{At} B V(t)\red{^2} dt.
\]
\red{Since $V(t)^2=V(t)$ for $v(t)\in\{0, 1\}^m$,}
Problem~\ref{prob:L0} is equivalent to the following problem:
\begin{equation}
\begin{aligned}
	& \underset{v}{\text{maximize}}
	& & J_1(v)\\
	& \text{subject to}
	& & v(t) \in \{0, 1\}^m \quad \forall t\in[0, T], \\
	& & & \|v_j\|_{L^0} \leq\alpha_j \quad \forall j \in\{1,2,\dots,m\},\\
	& & & \|v(t)\|_{\ell^0} \leq\beta \quad \forall t\in[0, T],
\end{aligned}
\label{prob:L0-2}
\end{equation}
where 
\begin{align}
	&J_1(v) \triangleq \int_{0}^{T} \bigl[ f_1(t), f_2(t), \dots, f_m(t)\bigr] v(t)~dt,\notag\\
	&f_j(t)\triangleq b_{j}^{\top}e^{A^{\top}t} e^{At}b_j, \quad j=1,2,\dots,m,\label{eq:f_j}
\end{align}
and $b_j$ is the $j$-th column of $B$.
\end{rmk}

\begin{rmk}
Compared to existing works~\cite{ITKKLCSS18,Ols19},
our formulation considers the component-wise $L^0$ norm $\|v_j\|_{L^0}$
and includes the $\ell^0$ norm of inputs,
by which we can adjust each $L^0$ cost of control variables and the number of control nodes.
While the optimal solution in the previous works is shown to be constructed 
from the top slice of the functions $f_1(t), f_2(t), \dots, f_m(t)$ defined by \eqref{eq:f_j} 
by using a rearrangement~\cite{Ols19},
this property does not hold for our optimization problem.
This will be illustrated in the example section.
\end{rmk}

\section{Analysis}
\label{sec:analysis}

The convex relaxation problem of Problem~\ref{prob:L0} is defined as follows:

\begin{problem}
\label{prob:L1}
Given $A\in\mathbb{R}^{n\times n}$, $B\in\mathbb{R}^{n\times m}$, $T>0$, $\beta\in\{1,2,\cdots,m-1\}$, and $\alpha_j\in(0, T]$, $j=1,2,\dots,m$,
find a function 
$v(t)\triangleq[v_1(t), v_2(t), \dots, v_m(t)]^{\top}$
that solves
\begin{equation*}
\begin{aligned}
	& \underset{v}{\text{maximize}}
	& & J_1(v) \\  & \text{subject to}
	& & v(t) \in [0, 1]^m \quad \forall t\in[0, T], \\
	& & & \|v_j\|_{L^1} \leq\alpha_j \quad \forall j \in\{1,2,\dots,m\},\\
	& & & \|v(t)\|_{\ell^1} \leq\beta \quad \forall t\in[0, T].
\end{aligned}
\end{equation*}
\end{problem}

The set of all functions that satisfy the constraints of an optimization problem is called {\em feasible set}.
Let us denote the feasible set of Problem~\ref{prob:L0} and Problem~\ref{prob:L1} by $\mathcal{V}_0$ and $\mathcal{V}_1$, i.e.,
\begin{align*}
	&\mathcal{V}_0 \triangleq
		\{v:~v(t)\in\{0,1\}^m ~ \forall t,~
	 	\|v_j\|_{L^0} \leq\alpha_j ~ \forall j,~
	  	\|v(t)\|_{\ell^0}\leq \beta ~ \forall t\},\\
	&\mathcal{V}_1 \triangleq
		\{v:~v(t)\in[0,1]^m ~ \forall t,~
		\|v_j\|_{L^1} \leq\alpha_j ~ \forall j,~
		\|v(t)\|_{\ell^1}\leq \beta ~ \forall t\}.
\end{align*}
Note that $\mathcal{V}_0\subset\mathcal{V}_1$, 
since $\|v_j\|_{L^1}=\|v_j\|_{L^0}$ for all $j$ and $\|v(t)\|_{\ell^1} = \|v(t)\|_{\ell^0}$ on $[0, T]$ 
for any measurable function $v$ satisfying $v(t)\in\{0,1\}^m$ for all $t$.
Then, we first show the discreteness of solutions of Problem~\ref{prob:L1},
which guarantees that the optimal solutions of Problem~\ref{prob:L1} belong to the set $\mathcal{V}_0$.

\begin{theorem}[discreteness]
\label{thm:discrete-L1}
Define functions $f_j$ by \eqref{eq:f_j}. 
If $f_j$ and $f_j -f_i$ are not constant on $[0, T]$ for all $i, j\in\{1,2,\dots,m\}$ with $j \neq i$,
then the solution of Problem~\ref{prob:L1} is unique\footnote{On any set of measure zero, any solution can take any values without loss of optimality. Hence, throughout the paper, when we say the optimal solution is {\em unique}, we mean up to sets of measure zero.} 
and it takes only the values in the binary set $\{0, 1\}$ almost everywhere.
\end{theorem}
\begin{IEEEproof}
Note that, for any $v$ such that $v(t)\in [0, 1]^m$ on $[0, T]$, we have
\begin{align*}
	&\|v(t)\|_{\ell^1}
	=  \sum_{j=1}^{m} |v_j(t)|
	= \sum_{j=1}^{m} v_j(t)
	= [1, 1, \dots, 1] v(t),\\
	&\|v_j\|_{L^1}
	=  \int_{0}^{T} |v_j(t)|dt
	= \int_{0}^{T} v_j(t)dt.
\end{align*}
Then, for each $j$, the value $\|v_j\|_{L^1}$ is equal to the final state $y_j(T)$ of the system
$\dot{y_j}(t) = v_j(t)$ with $y_j(0)=0$.
Hence, Problem~\ref{prob:L1} is equivalently expressed as follows:
\begin{equation}
\begin{aligned}
	\begin{split}
	& \underset{v}{\text{maximize}}
	& & J_1(v) \\
	& \text{subject to}
	& &   \dot{y}(t) = v(t), \quad y(0)=0,\\
	& & & y_j(T)\leq\alpha_j \quad \forall j\in\{1,2,\dots,m\},\\
	& & & v(t) \in [0, 1]^m \quad \forall t\in[0, T],\\
	& & & [1, 1, \dots, 1] v(t) \leq \beta \quad \forall t\in[0, T].
	\end{split}
\end{aligned}
\label{prob:L1-2}
\end{equation}
This is an optimal control problem to which 
Pontryagin's maximum principle~\cite[Theorem 22.2]{Cla13} is applicable.

Let the process $(y^{\ast}, v^{\ast})$ be a local maximizer of the problem~\eqref{prob:L1-2}.
Then, it follows from the maximum principle that 
there exists a constant $\eta$ equal to $0$ or $1$ and \red{an arc} $q:[0, T]\to{\mathbb{R}}^{m}$ satisfying the following conditions:
\begin{enumerate}
\item[(i)] the nontriviality condition:
\begin{equation}
	(\eta, q(t)) \neq 0 \quad \forall t\in[0, T], 
\label{eq:nontrivial_cond}
\end{equation}
\item[(ii)] the transversality condition:
\begin{equation}
	-q(T) \in N_{S}^{L} (y^\ast(T)), 
\label{eq:trans_cond}
\end{equation}
where $S \triangleq \{a\in\mathbb{R}^m: a_j \leq \alpha_j, ~ \forall j\}$,
\item[(iii)] the adjoint equation for almost every $t\in[0, T]$:
\begin{equation}
	-{\dot q}(t) = D_y H^{\eta}(t, y^{\ast}(t), q(t), v^{\ast}(t)), 
\label{eq:adjoint_cond}
\end{equation}
where
$D_y H^\eta$ is the derivative of the function $H^\eta$ at the second variable $y$, 
and $H^{\eta}: [0, T] \times \mathbb{R}^m \times \mathbb{R}^m \times \mathbb{R}^m$
is the Hamiltonian function associated to the problem~\eqref{prob:L1-2}, which is defined by
\[
	H^{\eta}(t, y, q, v) \triangleq q^\top v + \eta [f_1(t), f_2(t), \dots, f_m(t)] v,
\]
\item[(iv)] the maximum condition for almost every $t\in[0, T]$:
\begin{equation}
	H^{\eta}(t, y^{\ast}(t), q(t), v^{\ast}(t)) =  \sup_{v\in\mathbb{V}} H^{\eta}(t, y^{\ast}(t), q(t), v),
\label{eq:max_cond}
\end{equation}
where $\mathbb{V}\triangleq\{v\in[0, 1]^m: \sum_{j=1}^{m} v_j\leq\beta\}$.
\end{enumerate}

It follows from \eqref{eq:adjoint_cond} and \cite[Theorem~6.41]{Cla13}
that there exists a constant $\gamma\in\mathbb{R}^m$ such that $q(t)=\gamma$ on $[0,T]$,
since our Hamiltonian does not depend on the second variable $y$.
Then, from~\eqref{eq:trans_cond},
there exist sequences $\{\xi_i\}\subset\mathbb{R}^m$ and $\{\omega_i\}\subset\mathbb{R}^m$ 
such that
\begin{align}
\begin{split}
	&-\gamma = \lim_{i\to\infty}\xi_i, \quad
	 \xi_i \in N_S^P(\omega_i) \quad \forall i \in \mathbb{N}, \\
	&\lim_{i\to\infty} \omega_i = y^\ast(T),\quad
	\omega_i \in S \quad \forall i \in \mathbb{N}.
\end{split}
\label{eq:xi_wi}
\end{align}
For any $i\in\mathbb{N}$, by definition, there exists a constant $\sigma_i\geq0$ such that 
\begin{equation}
	\xi_i^\top (w - \omega_i) \leq \sigma_i \|w - \omega_i\|^2 \quad \forall w \in S.
\label{eq:proximal_normal}
\end{equation}
Fix any $\varepsilon>0$,  $j_0\in\{1, 2, \dots, m\}$, and $i\in\mathbb{N}$.
Take $w \in S$ such that 
$w_{j_0} = \omega_{i}^{(j_0)} - \varepsilon$ and $w_j = \omega_{i}^{(j)}$ for $j\neq j_{0}$,
where $\omega_{i}^{(j)}$ denotes the $j$-th component of $\omega_{i}\in\mathbb{R}^m$.
From the inequality~\eqref{eq:proximal_normal}, we have 
$\xi_i^{(j_0)} \geq - \sigma_i \varepsilon$, 
which gives $\xi_i^{(j)}\geq0$ for all $j\in\{1, 2, \dots, m\}$ and $i\in\mathbb{N}$
from the arbitrariness of $\varepsilon>0$, $j_0$, and $i$.
Hence, $\gamma_j \leq 0$ for all $j$ by \eqref{eq:xi_wi}.

In addition, we have $\gamma_j (y_j^{\ast}(T) - \alpha_j) = 0$ for all $j$.
Indeed, if $y_{j_0}^{\ast}(T) - \alpha_{j_0} \neq 0$ for some $j_0$, i.e., $y_{j_0}^{\ast}(T) < \alpha_{j_0}$,
then there exists $N\in\mathbb{N}$ such that 
$\omega_i^{(j_0)} < \alpha_{j_0}$ for all $i \geq N$ from \eqref{eq:xi_wi}.
Hence, $w \in \mathbb{R}^m$ defined by 
$w_{j_0} = \omega_i^{(j_0)}+\varepsilon$ and $w_j =\omega_i^{(j)}$ for $j\neq j_0$
satisfies $w\in S$ for sufficiently small $\varepsilon>0$.
Then, we have $\xi_i^{(j_0)} \leq \sigma_i \varepsilon$ from \eqref{eq:proximal_normal}.
This with $\xi_i^{(j_0)}\geq0$ implies $\xi_i^{(j_0)}=0$, which holds for all $i>N$. 
From \eqref{eq:xi_wi}, we have $\gamma_{j_0}=0$,
and thus $\gamma_{j_0}(y_{j_0}^{\ast}(T) - \alpha_{j_0}) = 0$.
Finally, the supremum in \eqref{eq:max_cond} is attained by a point in $\mathbb{V}$, 
since the right hand side is a linear function of $v$ and  $\mathbb{V}$ is a closed set.

In summary, the necessary conditions are given by 
\begin{align}
	&(\eta, \gamma) \neq 0, \label{eq:nonzero}\\ 
	&\gamma_j \leq 0 \quad \forall j\in\{1,2,\dots,m\}, \label{eq:nonpositive}\\
	&\gamma_j \left(y_j^{\ast}(T) - \alpha_j\right) = 0 \quad \forall j\in\{1,2,\dots,m\}, \label{eq:slack}\\
	&v^{\ast}(t) \in \argmax_{v\in \mathbb{V}} \sum_{j=1}^{m} \left(\eta f_{j}(t) + \gamma_j \right) v_{j}\label{eq:max}
\end{align}
almost everywhere.
We here claim that $\eta=1$, which can be observed as follows:
Assume $\eta=0$.
From~\eqref{eq:nonzero} and \eqref{eq:nonpositive}, there exists $j$ such that $\gamma_j < 0$.
Then, $v_j^\ast(t)=0$ almost everywhere by \eqref{eq:max}.
This implies $y_j^\ast(T)=0$ by the dynamics $\dot{y}_j = v_j$ with the initial condition $y_j(0)=0$.
Then, we have $\gamma_j (y_j^{\ast}(T) - \alpha_j) =-\gamma_j \alpha_j\neq0$,
which contradicts to \eqref{eq:slack}.
Thus, $\eta=1$.

Note that $\mu_L(E_{j})=0$ for all $j$, 
where $E_j \triangleq \{t\in[0, T]: f_{j}(t) + \gamma_{j} = 0\}$.
Indeed, if $\mu_{L}(E_{j})>0$ for some $j$, 
then an analytic function $\phi_{j} (t) \triangleq f_{j}(t) + \gamma_j$ takes zero on $E_{j}$ with a positive measure, 
which implies $\phi_{j}\equiv 0$ by \cite[Chapter~1]{Kra02} and contradicts to the assumption that $f_{j}$ is not constant.
Note also that 
there exist some functions $j_k(t): [0,T]\to\{1,2,\dots,m\}$, $k=1,2,\dots,m$,
such that we have
$\{j_1(t), j_2(t),\dots, j_m(t)\} = \{1, 2,\dots, m\}$ and 
\[
	f_{j_1(t)}(t) + \gamma_{j_1(t)} > \cdots > f_{j_m(t)}(t) + \gamma_{j_m(t)}
\]
almost everywhere,
which follows from the assumption that 
analytic functions $f_j - f_i$ are not constant for all $j\neq i$.
Hence, for almost every $t$,
\[
	v_j^\ast(t)=
	\begin{cases}
	1, &\mbox{if~} j \in \Lambda_\beta(t)\cap\Omega(t), \\
	0, &\mbox{otherwise}
	\end{cases}
\]
and the optimal solution is unique, where 
\begin{align*}
	&\Lambda_\beta(t)\triangleq \{j_1(t), j_2(t), \dots, j_\beta(t)\},\\
	&\Omega(t) \triangleq \{k\in\{1,2,\dots,m\}: f_{j_k(t)}(t) + \gamma_{j_k(t)}>0\}.
\end{align*}
This completes the proof.
\end{IEEEproof}


The following theorem is the main result, 
which shows the equivalence between Problem~\ref{prob:L0} and Problem~\ref{prob:L1}.

\begin{theorem}[equivalence]
\label{thm:L0-L1}
Define functions $f_j$ by \eqref{eq:f_j}, 
and assume $f_j$ and $f_j -f_i$ are not constant on $[0, T]$ for all $i, j\in\{1,2,\dots,m\}$ with $j \neq i$.
Denote the set of all solutions of Problem~\ref{prob:L0} and Problem~\ref{prob:L1} by 
${\mathcal{V}}_{0}^{\ast}$ and ${\mathcal{V}}_{1}^{\ast}$, respectively.
If the set ${\mathcal{V}}_{1}^{\ast}$ is not empty,
then ${\mathcal{V}}_{0}^{\ast}={\mathcal{V}}_{1}^{\ast}$ holds\footnote{Precisely, the equality means in the sense of equality of sets of equivalence classes, since any two functions that are equal almost everywhere are identified.}.
\end{theorem}
\begin{IEEEproof}
It follows from Theorem~\ref{thm:discrete-L1} that 
the solution of Problem~\ref{prob:L1} is unique and
it takes only the values $0$ and $1$ almost everywhere.
Let $\hat{v}\triangleq [\hat{v}_1, \hat{v}_2,\dots, \hat{v}_m]^{\top}\in\mathcal{V}_1^\ast$.
Note that the null set
$\cup_{j=1}^{m} \{t\in[0, T]: \hat{v}_j(t)\not\in\{0, 1\}\}$
does not affect the cost,
and hence we can adjust the variables so that $\hat{v}_j(t) \in \{0, 1\}$ on $[0, T]$ for all $j$,
without loss of the optimality.
We have
\begin{align*}
	\|\hat{v}(t)\|_{\ell^1} = \|\hat{v}(t)\|_{\ell^0},\quad
	\|\hat{v}_j\|_{L^1} = \|\hat{v}_j\|_{L^0}
\end{align*}
for all $j$, where we used the discreteness of $\hat{v}$.
Since $\hat{v}\in\mathcal{V}_1$, 
we have $\|\hat{v}(t)\|_{\ell^0}\leq\beta$ and $\|\hat{v}_j\|_{L^0}\leq\alpha_j$ for all $t$ and $j$.
Thus, $\hat{v}\in\mathcal{V}_0$.
This with $\mathcal{V}_0 \subset \mathcal{V}_1$ and $\hat{v}\in\mathcal{V}_1^\ast$
immediately gives
\[
	J_1(\hat{v}) 
	\leq \max_{v\in \mathcal{V}_0} J_1(v) 
	\leq \max_{v\in \mathcal{V}_1} J_1(v) 
	= J_1(\hat{v}).
\]
Hence, we have 
\begin{equation}
	J_1(\hat{v}) = \max_{v\in \mathcal{V}_0} J_1(v),
\label{eq:L0-L1}
\end{equation}
which implies $\hat{v}\in\mathcal{V}_0^\ast$.
Hence, 
${\mathcal{V}}_{1}^{\ast}\subset{\mathcal{V}}_{0}^{\ast}$ 
and ${\mathcal{V}}_{0}^{\ast}$ is not empty.

Next, take any 
$\tilde{v} \in{\mathcal{V}}_{0}^{\ast}$.
Note that $\tilde{v} \in \mathcal{V}_1$, 
since $\mathcal{V}_0^\ast \subset \mathcal{V}_0 \subset \mathcal{V}_1$.
In addition, it follows from \eqref{eq:L0-L1} that
$J_1(\tilde{v}) = J_1(\hat{v})$.
Therefore, $\tilde{v}\in{\mathcal{V}}_{1}^{\ast}$,
which implies ${\mathcal{V}}_{0}^{\ast}\subset{\mathcal{V}}_{1}^{\ast}$.
This gives ${\mathcal{V}}_{0}^{\ast}={\mathcal{V}}_{1}^{\ast}$.
\end{IEEEproof}

The existence of optimal solutions of Problem~\ref{prob:L1} is assumed 
in Theorem~\ref{thm:discrete-L1} and Theorem~\ref{thm:L0-L1}.
Although we omit the proof due to the page limitation,
we can show the existence in a similar way to \cite[Lemma~1]{ITKKTAC18}.

\begin{conjecture}[existence] \label{thm:existence}
For any $A\in\mathbb{R}^{n \times n}$, $B\in\mathbb{R}^{n \times m}$, $T>0$, 
$\beta\in\{1,2,\dots,m-1\}$, and $\alpha_{j}\in(0, T]$, $j=1,2,\dots,m$,
optimal solutions of Problem~\ref{prob:L1} exist.
\end{conjecture}

\subsection{Example}
\label{subsec:example}

\begin{figure}[t]
  \centering
    \includegraphics[width=0.8\linewidth]{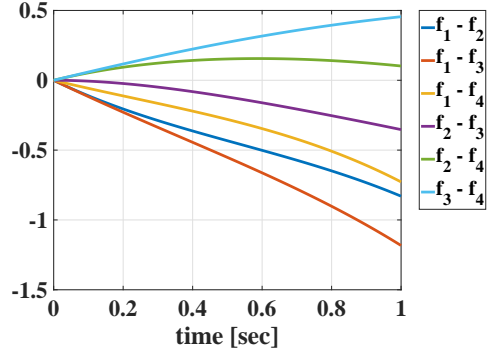}
  \caption{Functions $f_i(t) - f_j(t)$ for all $i,j$ with $i\neq j$}
\label{fig:difference}
\end{figure}

\begin{figure}[t]
  \centering
    \includegraphics[width=0.7\linewidth]{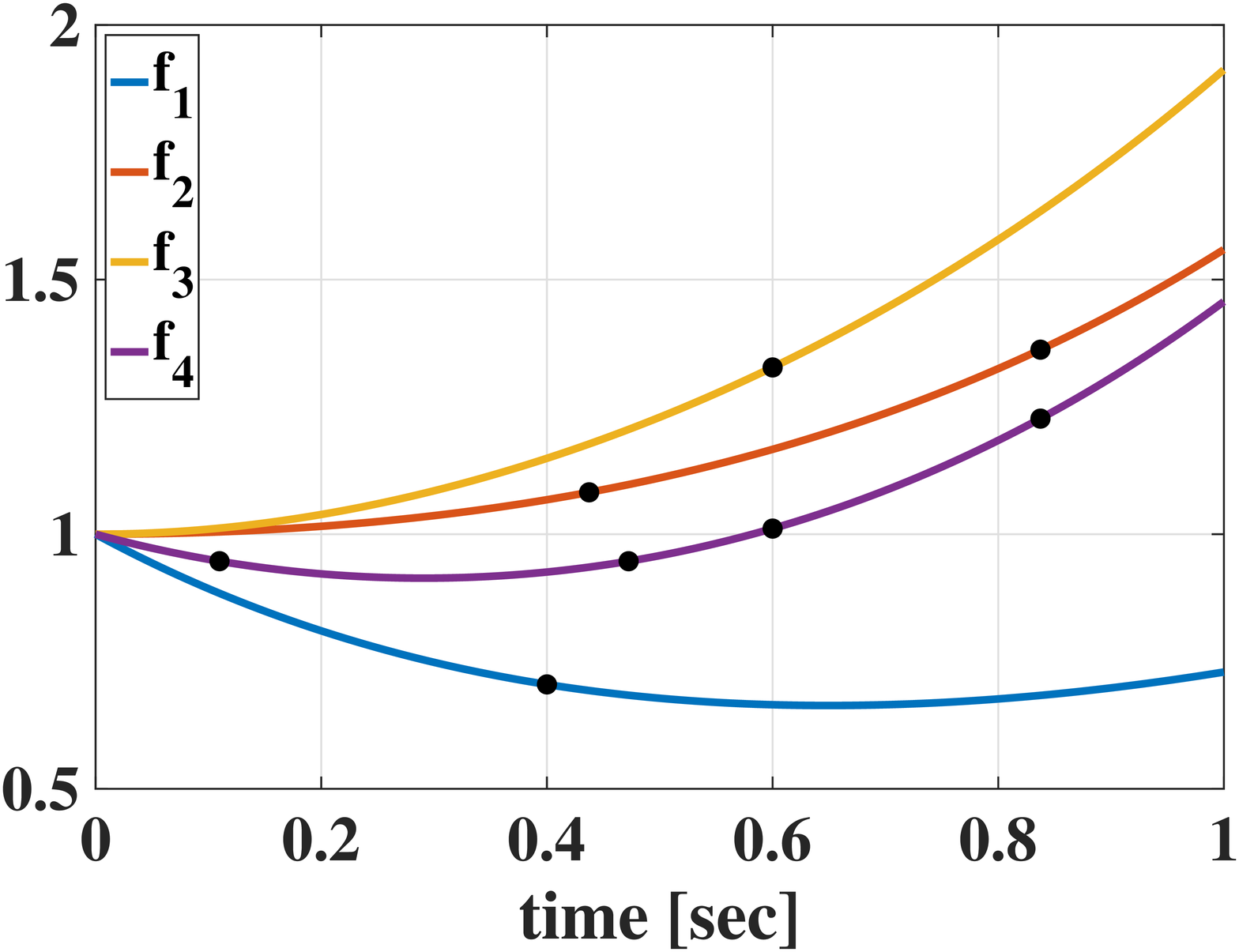}
    \includegraphics[width=0.7\linewidth]{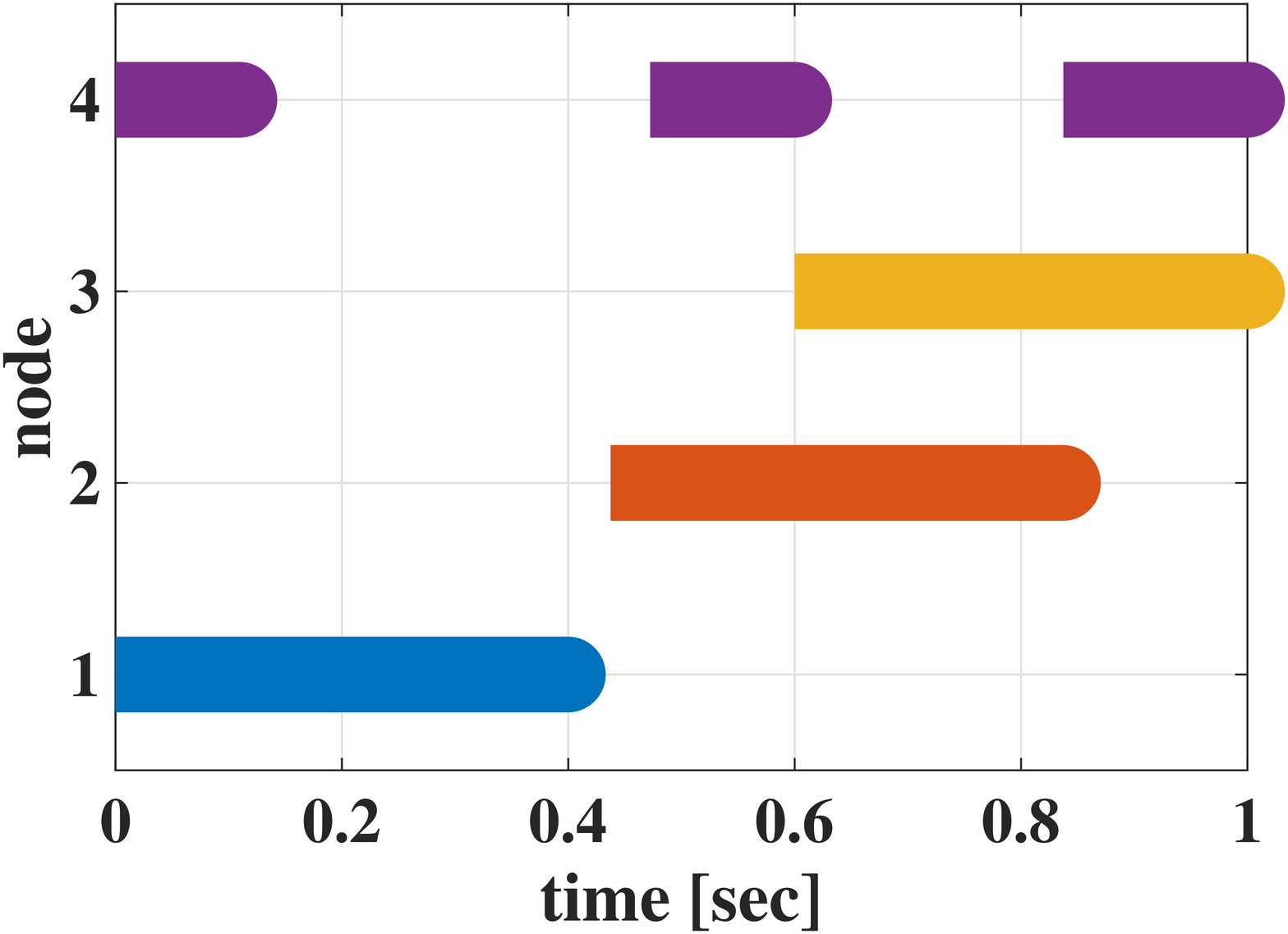}
  \caption{Functions $f_j(t)$ with switching points (top) and the proposed control node scheduling (bottom)}
\label{fig:node_ex1}
\end{figure}

\begin{figure}[htb]
  \centering
    \includegraphics[width=0.7\linewidth]{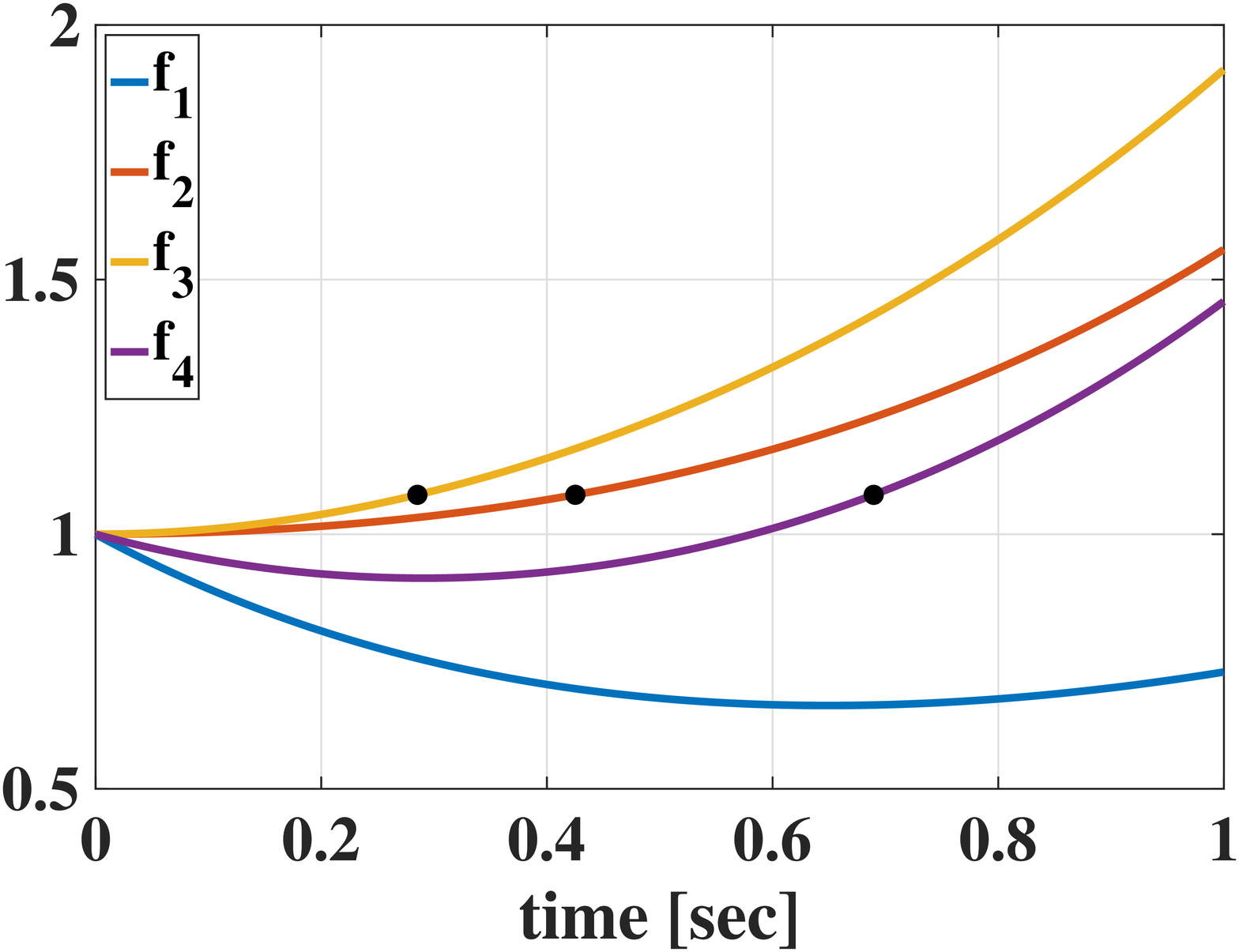}
    \includegraphics[width=0.7\linewidth]{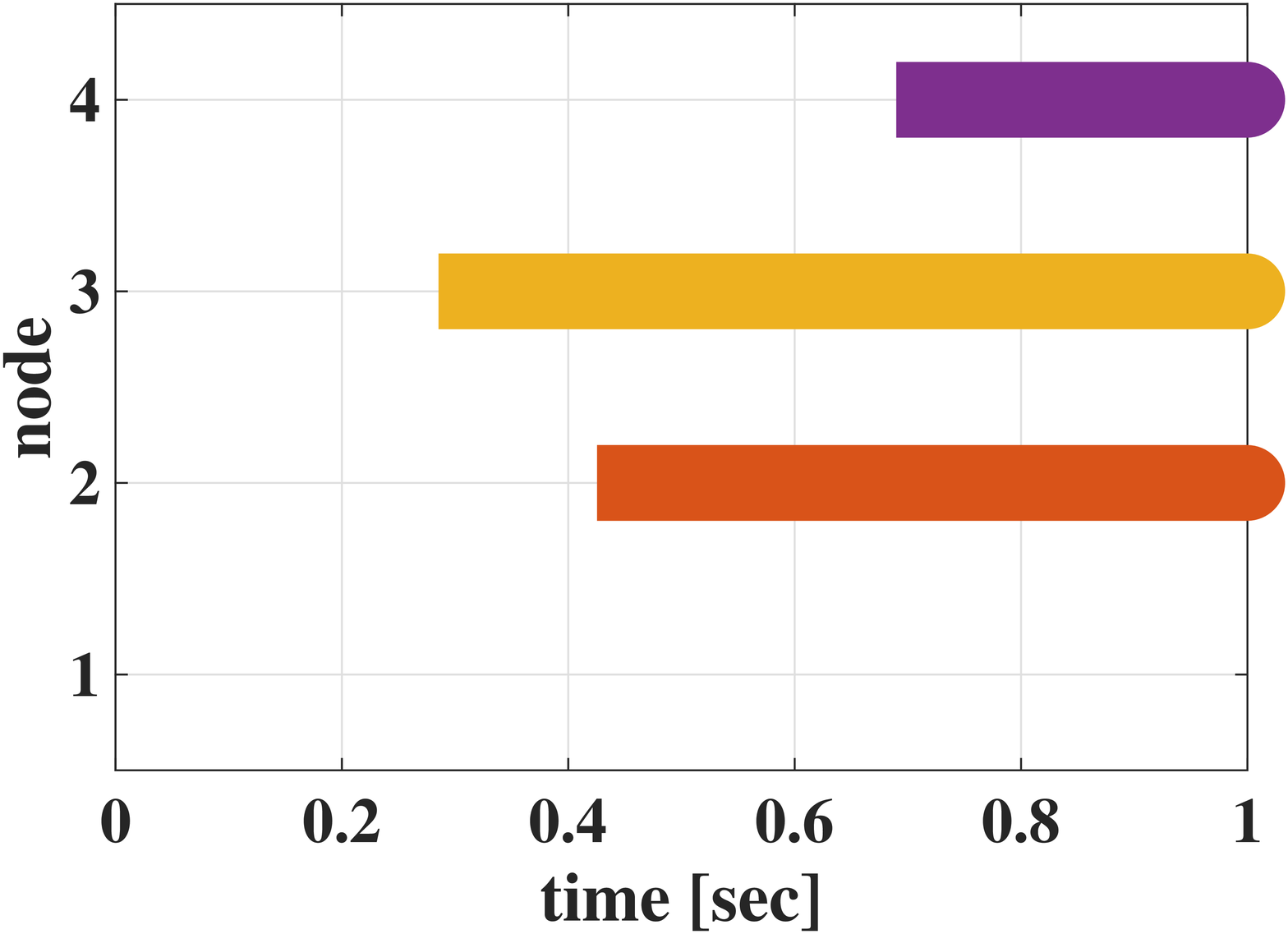}
  \caption{Functions $f_j(t)$ with switching points (top) and the control node scheduling in~\cite{ITKKLCSS18,Ols19} (bottom)}
\label{fig:node_ex2}
\end{figure}

This section gives an example of our node scheduling.
We consider a network model~\eqref{eq:system2} consisting of $4$ nodes with
\[
A =
\begin{bmatrix}
	-0.6 & 0 & -0.6 & 0.2\\
	-0.5 & 0 & 0 & 0.4\\
	1 & 0.6 & 0 & 0.5\\
	0 & 0 & 0.9 & -0.3
\end{bmatrix}
\]
and $B = I_{4}$, which is the identity matrix of dimension $4$.
For this network, we simulated our proposed method with $T=1$, $\alpha_j=0.4$, $j=1,2,3,4$, and $\beta=2$.
In this example, each node can become a control node since $B=I_{4}$, 
but the $L^0$ and $\ell^0$ constraints impose us to 
select at most 2 control nodes at each time and provide a control input to each node at most $0.4$ sec.
Note that the functions $f_i(t)-f_j(t)$ and $f_j(t)$ 
are shown in Fig.~\ref{fig:difference} and the upper panel of Fig.~\ref{fig:node_ex1},
by which we can confirm the equivalence between Problem~\ref{prob:L0} and Problem~\ref{prob:L1} from Theorem~\ref{thm:L0-L1}.
Then, we applied CVX \cite{cvx} in MATLAB, which is a software for convex optimization, to Problem~\ref{prob:L1}.

Fig.~\ref{fig:node_ex1} shows the resulting time series of control nodes on $[0, T]$.
Certainly, we can see that the set of control nodes depends on the time and satisfies both the $L^0$ and $\ell^0$ constraints.
Thus, we can find a finite number of essential nodes at each time and an essential time interval to provide control inputs.
For comparison, we also simulated the previous node scheduling proposed in \cite{ITKKLCSS18,Ols19},
which try to find a function $v(t)\in\mathbb{R}^m$ that solves
\begin{equation}
\begin{aligned}
	& \underset{v}{\text{maximize}}
	& & J(v) \\  & \text{subject to}
	& & v(t) \in \{0, 1\}^m \quad \forall t\in[0, T], \\
	& & & \|v\|_{L^0} \leq\alpha,
\end{aligned}
\label{prob:previous}
\end{equation}
where 
$\alpha\in(0, mT)$ is given.
Fig.~\ref{fig:node_ex2} shows the optimal scheduling when $\alpha=1.6$, which is equal to the value $\sum_{j=1}^{4}\alpha_j$.
Compared to our framework,
the problem~\eqref{prob:previous} does not include the $\ell^0$ and the component-wise $L^0$ constraints.
Indeed, the optimal solution selects more than 2 control nodes on an interval and tends to select a particular node.
Note also that the feasible set of the problem~\eqref{prob:previous} includes the set $\mathcal{V}_0$ 
due to the absence of the $\ell^0$ constraint, 
by which the optimal value is greater than that of Problem~\ref{prob:L0},
where the optimal values of problem~\eqref{prob:previous} and Problem~\ref{prob:L0} 
are $2.1287$ and $1.8957$, respectively.
Finally, the values $f_j(t)$ at the switching instances of control nodes are all equal,
i.e., we have to take the ``top slice" of the functions $f_j$.
Actually, this property is shown in~\cite{Ols19}.
On the other hand, 
our optimal solution can not be obtained by the simple method, as shown in Fig.~\ref{fig:node_ex1}.
Indeed, node $4$ is chosen as a control node around the final time $T$,
although the value $f_4(t)$ is not ranked in the top $2$ among $f_1(t), \dots, f_4(t)$ at the time.


\section{Application to Rebalancing of Mobility Networks}
\label{sec:application_to_rebalancing}

\subsection{Problem Formulation}

In this section, we consider a rebalancing problem on the mobility network of a sharing system with one-way trips. This system is first modeled, 
and the rebalancing problem is then formulated by the main problem tackled in this paper.

The structure of the network is modeled, based on 
\cite{Calafiore2019}, as follows.
Let $\mathcal{S} = \{1,2,\ldots,s\}$ be the set of the indexes of $s \in \mathbb{N}$ stations, where vehicles are parked.
\red{Assume that customers take the service according to the Poisson process at a rate $g_{ij}(t) \ge \mathbb{R}$,
the number of the demand for travels from station $j \in \mathcal{S}$ to $i \in \mathcal{S}$ per time.
In addition, assume that the time it takes to travel from station $j$ to $i$ follows the exponential distribution with an average $\tau_{ij}(t) > 0$.
Let $u_{ij}(t) \ge 0$ be 
the rate at which the vehicles are rebalanced by the staff
from station $j$ to $i$.
Let $f_{ij}(t) \ge 0$ be the expectation of the number of vehicles traveling from station $j$ to $i$,
and it varies according to
}
\begin{equation}
    \dot f_{ij}(t) = -\gamma_{ij}(t)f_{ij}(t) + g_{ij}(t) + u_{ij}(t)
 \label{eq_sk:flow_dynamics}
 ,
\end{equation}
\red{where $\gamma_{ij}(t) = 1/\tau_{ij}(t)$ corresponds to the arrival rate of vehicles at station $i$ from $j$ per time.
Let $v_i(t) \ge 0$ denote the expectation of the number of vehicles parked at station $i$,
and it varies according to
\begin{equation}
 \dot v_i(t) =  \sum_{j=1}^s \gamma_{ij}(t) f_{ij}(t)-\sum_{j=1}^s (g_{ji}(t) + u_{ji}(t))
 \label{eq_sk:dot_v_i}
.
\end{equation}
}
\blue{See \cite{arXivSM} for detailed derivation of the system.}

Dynamic pricing is applied to this system, which is modeled as follows.
Let $p_{ij}(t) \in \mathbb{R}$ be the price for the rent of a vehicle from station $j$ to $i$.
According to an economic model \cite{Parkin17},
the amount of the vehicles in service depends on the price as follows:
\begin{equation}
    g_{ij}(t) = \bar g_{ij}(t) - \theta_{ij}(t) p_{ij}(t)
    \label{eq_sk:e_ij}
    ,
\end{equation}
where $\bar g_{ij}(t) \red{\ge 0}$ is the expected demand of the vehicles when $p_{ij}(t) = 0$, and $\theta_{ij}(t) \red{> 0}$ denotes the price elasticity.
Assume that the price $p_{ij}(t)$ is determined by the amount of the vehicles at stations $i$ and $j$
according to the following rule:
\begin{equation}
    p_{ij}(t) = \bar p_{ij}(t) + \lambda_{ij}(t)v_i(t)-\lambda_{ji}(t)v_j(t)
  \label{eq_sk:p_ij}
  ,
\end{equation}
where $\bar p_{ij}(t) \red{\ge 0}$ is a standard price,
and $\lambda_{ij}(t) \ge 0$ and  $\lambda_{ji}(t) \ge 0$
are the coefficients to adjust the price according to the amount of vehicles at stations $i$ and $j$, respectively.
Under this rule, as the amount $v_i(t)$ ($v_j(t)$) of vehicles at station $i$ ($j$) becomes larger (smaller),
the price $p_{ij}(t)$ becomes higher to reduce the amount of vehicles entering station $i$ (leaving station $j$).
Let the standard price $\bar p_{ij}(t)$ be determined 
according to the expected demand $\bar g_{ij}(t)$ of vehicles
as follows:
\begin{equation}
    \bar p_{ij}(t) = \frac{\bar g_{ij}(t)}{\theta_{ij}(t)}
    \label{eq_sk:bar_p_bar_g}
    .
\end{equation}

From (\ref{eq_sk:e_ij}), (\ref{eq_sk:p_ij}), and (\ref{eq_sk:bar_p_bar_g}),
(\ref{eq_sk:flow_dynamics}) \red{and (\ref{eq_sk:dot_v_i})}
are reduced to
\begin{align}
  \dot f_{ij}(t) &= -\gamma_{ij}(t) f_{ij}(t) -\theta_{ij}(t)(\lambda_{ij}(t)v_i(t)-\lambda_{ji}(t)v_j(t)) 
  \notag\\
  & \qquad + u_{ij}(t)
 \label{eq_sk:flow_dynamics2}\\
 \red{\dot v_i(t)} & \red{ =  \sum_{j=1}^s 
  (\gamma_{ij}(t) f_{ij}(t) +\theta_{ji}(t)(\lambda_{ji}(t)v_j(t)-\lambda_{ij}(t)v_i(t))}
  \notag\\
  & \qquad \red{- u_{ji}(t))}
 \label{eq_sk:dot_v_i2}
,
\end{align}
respectively. 
By collecting (\ref{eq_sk:flow_dynamics2}) and \red{(\ref{eq_sk:dot_v_i2})} for all $i \in \mathcal{S}$,
we obtain the model of the mobility network as
\begin{equation}
  \dot x(t) = A(t)x(t) + B(t)u(t)
  \label{eq_sk:dot_x}.
\end{equation}
This is the system of dimension 
\red{$n = s^2$} 
with 
\red{$m = s^2-s$}  
inputs for
\begin{align}
x(t) &= [v_1(t)~v_2(t)~\cdots~v_s(t)~f_{12}(t)~f_{13}(t)~\cdots~\red{f_{s,s-1}(t)}]^\top,
\notag\\
u(t) &= [u_{12}(t)~u_{13}(t)~\cdots~\red{u_{s,s-1}(t)}]^\top,
\notag\\
 A(t) &= 
  \begin{bmatrix}
   \red{\Xi \Lambda(t)} & \Gamma(t)\\
   -\Lambda(t) & -\Delta(t)
  \end{bmatrix},~
 B(t) = 
 \begin{bmatrix}
   \red{-\Xi} \\
   I_m
\end{bmatrix}
 ,
\label{eq_sk:AB}
\end{align}
where
\begin{align*}
  \red{\Xi} & \red{= [\hat E_1~\hat E_2~\cdots~\hat E_s],~~\hat E_i = [e_1~\cdots~e_{i-1}~e_{i+1}~\cdots~e_s],}\\
 \Lambda(t) &= 
 \begin{bmatrix}
   \theta_{12}(t)(\lambda_{12}(t)e_1-\lambda_{21}(t)e_2)^\top\\
   \theta_{13}(t)(\lambda_{13}(t)e_1-\lambda_{31}(t)e_3)^\top\\
   \vdots\\
   \red{\theta_{s,s-1}(t)(\lambda_{s,s-1}(t)e_s-\lambda_{s-1,s}(t)e_{s-1})^\top}
 \end{bmatrix}
 \!,\\
 \Gamma(t) &= 
[e_1 \gamma_{12}(t) ~ e_1 \gamma_{13}(t) ~\cdots ~   \red{e_s\gamma_{s,s-1}(t)}]
 ,\\
  \Delta(t) & = \mathrm{diag}(\gamma_{12}(t),\gamma_{13}(t),\ldots,\red{\gamma_{s,s-1}(t)}),
\end{align*}
$I_m \in \mathbb{R}^{m \times m}$ is the identity matrix of dimension $m$,
and $e_i \in \mathbb{R}^s$ is the unit vector of dimension $s$ with the $i$th entry one.

For rebalancing,
the staff organizes $\beta \in \mathbb{N}$ teams to transport vehicles.
Each team goes to a station full of vehicles by their management car,
and some members of the team transfer vehicles to a vacant station along with the management car driven by other members. 
Then, the management car picks up the members to go to another station.
\red{Assume that the dynamics of staff is sufficiently faster than that of customers and that 
only one team can rebalance a vehicle on each route between stations
to distribute the staff around the area.}
Then, the number of rebalances at the same time is at most $\beta$, 
which is expressed as $\|u(t)\|_{\ell^{0}}\leq\beta$ for any $t$.
Without the loss of generality, the possible rebalance number is one
through standardization, 
which is expressed as $|u_{ij}(t)| \leq 1$ for any $i,j \in \mathcal{S}$ and any $t$, i.e.,
$\|u\|_{L^{\infty}} \leq 1$.
\red{Additionally, $u_{ij}(t) \ge 0$ has to be satisfied.}
\blue{See \cite{arXivSM} for dynamic models of rebalancing by the staff.}

Let $x_0 \in \mathbb{R}^n$ be the initial amounts of vehicles in the stations, which are unevenly distributed, 
and let $x_d \in \mathbb{R}^n$ be the desired terminal amounts to attain even distribution.
We design control input $u(t)$ to achieve the state $x(T) = x_d$ at the terminal time $T$ from the initial state $x(0) = x_0$
under the dynamics (\ref{eq_sk:dot_x}).
In summary, the balancing problem of the mobility network is formulated as follows.

\begin{problem}\label{prob:mobility}
Given $A(t)\in\mathbb{R}^{n\times n}$, 
$B(t)\in\mathbb{R}^{n\times m}$, 
$T>0$, 
$x_0\in\mathbb{R}^n$,
$x_d\in\mathbb{R}^n$,
and $\beta\in\{1,2,\cdots,m-1\}$, 
find a control $u$ that solves
\begin{equation*}
\begin{aligned}
	& \underset{u}{\text{minimize}}
	& & \|u\|_{L^0}\\
	& \text{subject to}
	& & \dot x(t) = A(t) x(t) + B(t) u(t) \quad \forall t\in[0, T], \\
	& & & x(0)=x_0, \quad x(T)=x_d, \\
	& & & \|u(t)\|_{\ell^0} \leq\beta \quad \forall t\in[0, T],\\
	& & & \red{u(t)\in[0, 1]^{m} \quad \forall t\in[0, T]}.
\end{aligned}
\end{equation*}
\end{problem}

We denote the feasible set of Problem~\ref{prob:mobility} by $\mathcal{U}_0$.
We also denote the state-transition matrix of $A(t)$ by $\Phi(t,\tau)$.
In other words, $\Phi(t,\tau)$ is the unique solution of the matrix differential equation
\begin{align*}
	\frac{d}{dt} \Phi(t,\tau) = A(t) \Phi(t,\tau),\quad
	\Phi(\tau,\tau) = I_{n},
\end{align*}
where $I_{n}$ is the identity matrix.

\begin{rmk}
\label{rmk:formulation_mobility}
Note that if $x_d=\Phi(T,0)x_0$, then the optimal control is trivial (i.e. zero control),
and hence we assume $x_d \neq \Phi(T,0)x_0$ throughout the paper.
Note also that we assume $\beta < m$, 
since if $\beta=m$ the optimal control is a standard time-sparse hands-off control 
discussed in~\cite{NagQueNes16,ITKKTAC18}.
\end{rmk}

\subsection{Analysis}
\label{sec:analysis_mobility}

We here introduce a convex relaxation problem of Problem~\ref{prob:mobility},
where the $L^0$ and $\ell^0$ norms are replaced by the $L^1$ and $\ell^1$ norms, respectively.

\begin{problem}\label{prob:convex}
Given $A(t)\in\mathbb{R}^{n\times n}$, 
$B(t)\in\mathbb{R}^{n\times m}$, 
$T>0$, 
$x_0\in\mathbb{R}^n$,
$x_d\in\mathbb{R}^n$,
and $\beta\in\{1,2,\cdots,m-1\}$,
find a control $u$ that solves
\begin{equation*}
\begin{aligned}
	& \underset{u}{\text{minimize}}
	& & \|u\|_{L^1}\\
	& \text{subject to}
	& & \dot x(t) = A(t) x(t) + B(t) u(t) \quad \forall t\in[0, T], \\
	& & & x(0)=x_0, \quad x(T)=x_d, \\
	& & & \|u(t)\|_{\ell^1} \leq\beta \quad \forall t\in[0, T],\\
	& & & \red{u(t)\in[0, 1]^{m} \quad \forall t\in[0, T]}.
\end{aligned}
\end{equation*}
\end{problem}

We denote the feasible set of Problem~\ref{prob:convex} by $\mathcal{U}_1$.
In other words,
\begin{align*}
	\mathcal{U}_0
	\triangleq \Bigg\{u:~& \Phi(T,0)x_0 + \int_{0}^{T} \Phi(T,t)B(t)u(t)dt = x_d, \\
	&\red{u(t)\in[0, 1]^{m}}, \|u(t)\|_{\ell^0}\leq \beta~\forall t\in[0, T] \Bigg\},\\
	\mathcal{U}_1
	\triangleq \Bigg\{u:~& \Phi(T,0)x_0 + \int_{0}^{T} \Phi(T,t)B(t)u(t)dt = x_d, \\
	&\red{u(t)\in[0, 1]^{m}}, \|u(t)\|_{\ell^1}\leq \beta~\forall t\in[0, T] \Bigg\}.
\end{align*}
Note that we have $\mathcal{U}_0\subset\mathcal{U}_1$, 
since $\|a\|_{\ell^1}\leq \|a\|_{\ell^0}$ for any \red{$a\in[0,1]^m$}.
Then, we first show the \red{discreteness} of the optimal solutions of Problem~\ref{prob:convex}.
The property guarantees that 
the optimal solutions of Problem~\ref{prob:convex} belong to the set $\mathcal{U}_0$,
which is illustrated in the proof of Theorem~\ref{thm:L0-L1_mobility}.
For this, we introduce an assumption on the impulse response $\Phi(\cdot,\cdot)B(\cdot)$ of the system.
\red{\begin{assumption}\label{ass:mobility}
For any nonzero $\rho\in\mathbb{R}^n$ and $\eta\in\{0, 1\}$, we have
\begin{equation*}
  \mu_L\left(\left\{t\in[0, T]: \rho^\top \Phi(T,t)b_j(t) = \eta \right\}\right) = 0 
\end{equation*}
for all $j$, and 
\begin{equation*}
  \mu_L\left(\left\{t\in[0, T]: \rho^\top \Phi(T,t)b_i(t) = \rho^\top \Phi(T,t)b_j(t) \right\}\right) = 0
\end{equation*}
for all $i, j$ with $i\neq j$,
where $b_j(t)\in\mathbb{R}^{n}$ is the $j$-th column vector in matrix $B(t)$.
\end{assumption}}

\begin{theorem}\label{thm:discrete-L1_mobility}
Under Assumption~\ref{ass:mobility},
the optimal solution of Problem~\ref{prob:convex} is unique 
and it takes only the values in the set $\red{\{0, 1\}}$ almost everywhere.
\end{theorem}
\begin{IEEEproof}
Let the process $(x^{\ast}, u^{\ast})$ be a local minimizer of Problem~\ref{prob:convex}.
Then, it follows from Pontryagin's maximum principle~\cite[Theorem 22.2]{Cla13} that 
there exists a constant $\eta$ equal to $0$ or $1$ and \red{an arc} $q:[0, T]\to{\mathbb{R}}^{n}$ satisfying the following conditions:
\begin{enumerate}
\item the nontriviality condition:
\begin{equation}
	(\eta, q(t)) \neq 0 \quad \forall t\in[0, T], 
\label{eq:nontrivial_cond}
\end{equation}
\item the adjoint equation for almost every $t\in[0, T]$:
\begin{equation}
	-{\dot q}(t) = D_x H^{\eta}(t, x^{\ast}(t), q(t), u^{\ast}(t)), 
\label{eq:adjoint_cond_mobility}
\end{equation}
where
$D_x H^\eta$ is the derivative of the function $H^\eta$ at the second variable $x$, 
and $H^{\eta}: [0, T] \times \mathbb{R}^n \times \mathbb{R}^n \times \mathbb{R}^m \to \mathbb{R}$  
is the Hamiltonian function associated to Problem~\ref{prob:convex}, 
which is defined by
\[
	H^{\eta}(t, x, q, u) \triangleq q^\top (A(t)x + B(t)u) - \eta \|u\|_{\ell^1},
\]
\item the maximum condition for almost every $t\in[0, T]$:
\begin{equation}
	H^{\eta}(t, x^{\ast}(t), q(t), u^{\ast}(t)) 
	=  \sup_{u\in\mathbb{U}} H^{\eta}(t, x^{\ast}(t), q(t), u),
\label{eq:max_cond_mobility}
\end{equation}
where 
$\red{\mathbb{U}\triangleq\{u\in[0, 1]^{m}: \sum_{j=1}^{m} u_{j} \leq \beta\}}$.
\end{enumerate}

Note that 
\[
	D_x H^{\eta}(t, x^{\ast}(t), q(t), u^{\ast}(t)) = A(t)^\top q(t).
\]
Hence, we have $q(t) = \Psi(t,T)q(T)$ on $[0,T]$ from \eqref{eq:adjoint_cond_mobility}, 
where $\Psi(\cdot,\cdot)$ is the state-transition matrix of $-A(t)^\top$.
Note also that the supremum in \eqref{eq:max_cond_mobility} is attained by a point in $\mathbb{U}$, 
since the right hand side is continuous on $u$ and the set $\mathbb{U}$ is closed.
Hence, we have 
\begin{align}
\begin{split}
	u^\ast(t) 
	& \in \argmax_{u\in\mathbb{U}} \left( q(t)^\top B(t) u -\eta \|u\|_{\ell^1} \right)\\
	& = \red{\argmax_{u\in\mathbb{U}} \sum_{j=1}^{m} \left( q(t)^\top b_j(t) - \eta \right) u_{j}}
\label{eq:necessary_optimal}
\end{split}
\end{align}
almost everywhere.

We here claim that $q(T)\neq0$.
Indeed, if $q(T) = 0$, then it follows from the nontriviality condition~\eqref{eq:nontrivial_cond} that $\eta=1$.
From \eqref{eq:necessary_optimal}, we have $u^\ast (t) = 0$ for almost all $t$ \blue{since $q(t)=0$ on $[0, T]$}. 
This implies $x_d=x^\ast(T)=\Phi(T,0)x_0$.
This contradicts the definition of $x_0$ and $x_d$ (see Remark~\ref{rmk:formulation_mobility}).
Hence, $q(T)\neq0$.
Since $\Psi(t,T)=\Phi(T,t)^\top$ for any $t$ from \cite{kai80},
we have
\begin{align}
	&\red{q(t)^\top b_j(t) \neq \eta \mbox{ for all } j,} \label{eq:discreteness_coff2} \\
	&\red{q(t)^\top b_i(t) \neq q(t)^\top b_j(t) \mbox{ for all } i, j \mbox{ with } i\neq j} \label{eq:discreteness_coff} 
\end{align}
for almost all $t\in[0,T]$ from Assumption~\ref{ass:mobility}.

In what follows, we show the \red{discreteness} property of $u^{\ast}(t)$ 
for both cases of $\eta=0$ and $\eta=1$.
For the characterization,
we here take functions $j_k(t): [0,T]\to\{1,2,\dots,m\}$, $k=1,2,\dots,m$, 
such that 
$\{j_1(t), j_2(t),\dots, j_m(t)\} = \{1, 2,\dots, m\}$ and 
\[
	\red{q(t)^\top b_{j_1(t)}(t) > q(t)^\top b_{j_2(t)}(t) > \cdots > q(t)^\top b_{j_m(t)}(t)}
\]
almost everywhere.
Note that the existence of $j_k(t)$ is guaranteed by~\eqref{eq:discreteness_coff}.
Define sets 
\begin{align*}
	&\Lambda_\beta(t) \triangleq \{j_1(t), j_2(t), \dots, j_\beta(t)\},\\
	&\red{\Omega_\eta(t) \triangleq \{k\in\{1, 2, \dots, m\}: q(t)^{\top} b_{j_{k}(t)}(t) > \eta\}}.
\end{align*}
\red{Then, it immediately follows that we have 
\[
	u_j^\ast(t)=
	\begin{cases}
	1, & \mbox{if~} j \in \Lambda_\beta(t)\cap \Omega_\eta(t), \\
	0, &\mbox{otherwise}
	\end{cases}
\]
almost everywhere, which completes the proof.}
\end{IEEEproof}

The following theorem is the main result, 
which shows the equivalence between Problem~\ref{prob:mobility} and Problem~\ref{prob:convex}.

\begin{theorem}
\label{thm:L0-L1_mobility}
Suppose Assumption~\ref{ass:mobility}.
Denote the set of all solutions of Problem~\ref{prob:mobility} and Problem~\ref{prob:convex} by 
${\mathcal{U}}_{0}^{\ast}$ and ${\mathcal{U}}_{1}^{\ast}$, respectively.
If the set ${\mathcal{U}}_{1}^{\ast}$ is not empty,
then ${\mathcal{U}}_{0}^{\ast}={\mathcal{U}}_{1}^{\ast}$.
\end{theorem}
\begin{IEEEproof}
Take any $\hat{u}\in\mathcal{U}_1^\ast$.
It follows from Theorem~\ref{thm:discrete-L1_mobility} that $\red{\hat{u}(t)\in\{0, 1\}^{m}}$ almost everywhere.
It follows from the \red{discreteness} of $\hat{u}$, we have
\begin{equation}
	\|\hat{u}(t)\|_{\ell^1}
	= \sum_{j=1}^m |\hat{u}_j(t)|
	= \sum_{\{j: \hat{u}_j(t) \neq 0\}}  1
	= \|\hat{u}(t)\|_{\ell^0}.
\label{eq:l1_optimal_number}
\end{equation}
Since $\hat{u}\in\mathcal{U}_{1}^\ast\subset\mathcal{U}_1$, 
we have $\|\hat{u}(t)\|_{\ell^0}\leq\beta$ by~\eqref{eq:l1_optimal_number}.
Thus, $\hat{u}\in\mathcal{U}_0$.
In addition, 
\begin{equation}
	\|\hat{u}\|_{L^1}
	= \int_{0}^{T} \|\hat{u}(t)\|_{\ell^1} dt
	= \int_{0}^{T} \|\hat{u}(t)\|_{\ell^0} dt
	= \|\hat{u}\|_{L^0}
\label{eq:L1_optimal}
\end{equation}
from~\eqref{eq:l1_optimal_number}.
Here, since $\|a\|_{\ell^1}\leq \|a\|_{\ell^0}$ for any $\red{a\in[0,1]^m}$,
we have $\mathcal{U}_0\subset\mathcal{U}_1$ and $\|u\|_{L^1}\leq \|u\|_{L^0}$ for any $u \in \mathcal{U}_0$.
Hence, for any $u \in \mathcal{U}_0$, we have
\begin{align*}
	\|\hat{u}\|_{L^0} 
	= \|\hat{u}\|_{L^1} 
	\leq \|u\|_{L^1} 
	\leq \|u\|_{L^0},
\end{align*}
where the first relation follows from~\eqref{eq:L1_optimal}
and the second relation follows from $u\in\mathcal{U}_1$ and the optimality of $\hat{u}$.
This implies $\hat{u}\in\mathcal{U}_0^\ast$.
Hence, $\mathcal{U}_1^\ast \subset \mathcal{U}_0^\ast$, and $\mathcal{U}_0^\ast$ is not empty.

We next take any $\tilde{u}\in\mathcal{U}_0^\ast$.
Note that $\tilde{u}\in\mathcal{U}_1$, since $\mathcal{U}_0^\ast\subset\mathcal{U}_0\subset\mathcal{U}_1$.
In addition,
\begin{align*}
	\|\tilde{u}\|_{L^1}
	\leq \|\tilde{u}\|_{L^0}
	\leq \|\hat{u}\|_{L^0}
	= \|\hat{u}\|_{L^1}
	\leq \|\tilde{u}\|_{L^1},
\end{align*}
where the first inequality follows from the fact that $\|a\|_{\ell^1} \leq \|a\|_{\ell^0}$ for any $\red{a\in[0, 1]^m}$,
the second inequality follows from $\hat{u}\in\mathcal{U}_0$ and the optimality of $\tilde{u}$,
the third equality follows from \eqref{eq:L1_optimal},
and the last inequality follows from $\tilde{u}\in\mathcal{U}_1$ and the optimality of $\hat{u}$.
This gives $\|\tilde{u}\|_{L^1} = \|\hat{u}\|_{L^1}$, 
which implies $\tilde{u}\in\mathcal{U}_1^\ast$.
Thus, we have $\mathcal{U}_0^\ast \subset \mathcal{U}_1^\ast$.
\end{IEEEproof}


We finally provide the existence of optimal controls of Problem~\ref{prob:convex} without the proof, 
which is confirmed in a similar way to \cite[Lemma~1]{ITKKTAC18}.
Since we have $\mathcal{U}_0 \subset \mathcal{U}_1$,
if the set $\mathcal{U}_0$ is not empty, 
then there exists an optimal control of Problem~\ref{prob:mobility} under Assumption~\ref{ass:mobility}.

\begin{conjecture}\label{thm:existence_mobility}
Given $(A(t), B(t), x_0, x_d, T, \beta)$,
the set $\mathcal{U}_1^\ast$ is not empty if and only if the set $\mathcal{U}_1$ is not empty,
where $\mathcal{U}_1^\ast$ is defined in Theorem~\ref{thm:L0-L1_mobility}.
\end{conjecture}

\red{\begin{rmk}\label{rmk:Linf}
Although we consider non-negative controls in Problem~\ref{prob:mobility},
with a slight modification we can show results similar to those presented above (i.e., Theorems~\ref{thm:discrete-L1_mobility} and \ref{thm:L0-L1_mobility}) for controls with an $L^{\infty}$ constraint $\|u\|_{L^{\infty}} \leq 1$.
To be more precise, under a suitable assumption,
the optimal control of a corresponding convex relaxation problem is unique 
and it takes only the values in $\{0, \pm 1\}$ almost everywhere, 
and this shows the equivalence between the original problem and the relaxation problem 
as illustrated in Theorem~\ref{thm:L0-L1_mobility}.
In this sense, our optimal control with multiple sparsity includes the existing sparse control in~\cite{NagQueNes16, ITKKTAC18} as a special case.
\end{rmk}}

\subsection{Numerical Study}
\begin{figure}[t]
  \centering
    \includegraphics[width=0.8\linewidth]{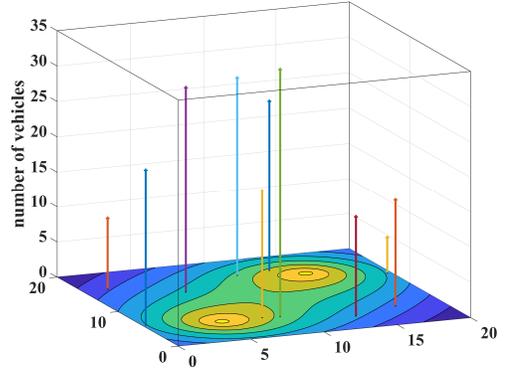}
  \caption{\red{Location of stations and the initial distribution}}
\label{fig:station_ex2}
\end{figure}

\begin{figure}[t]
  \centering
    \includegraphics[width=0.75\linewidth]{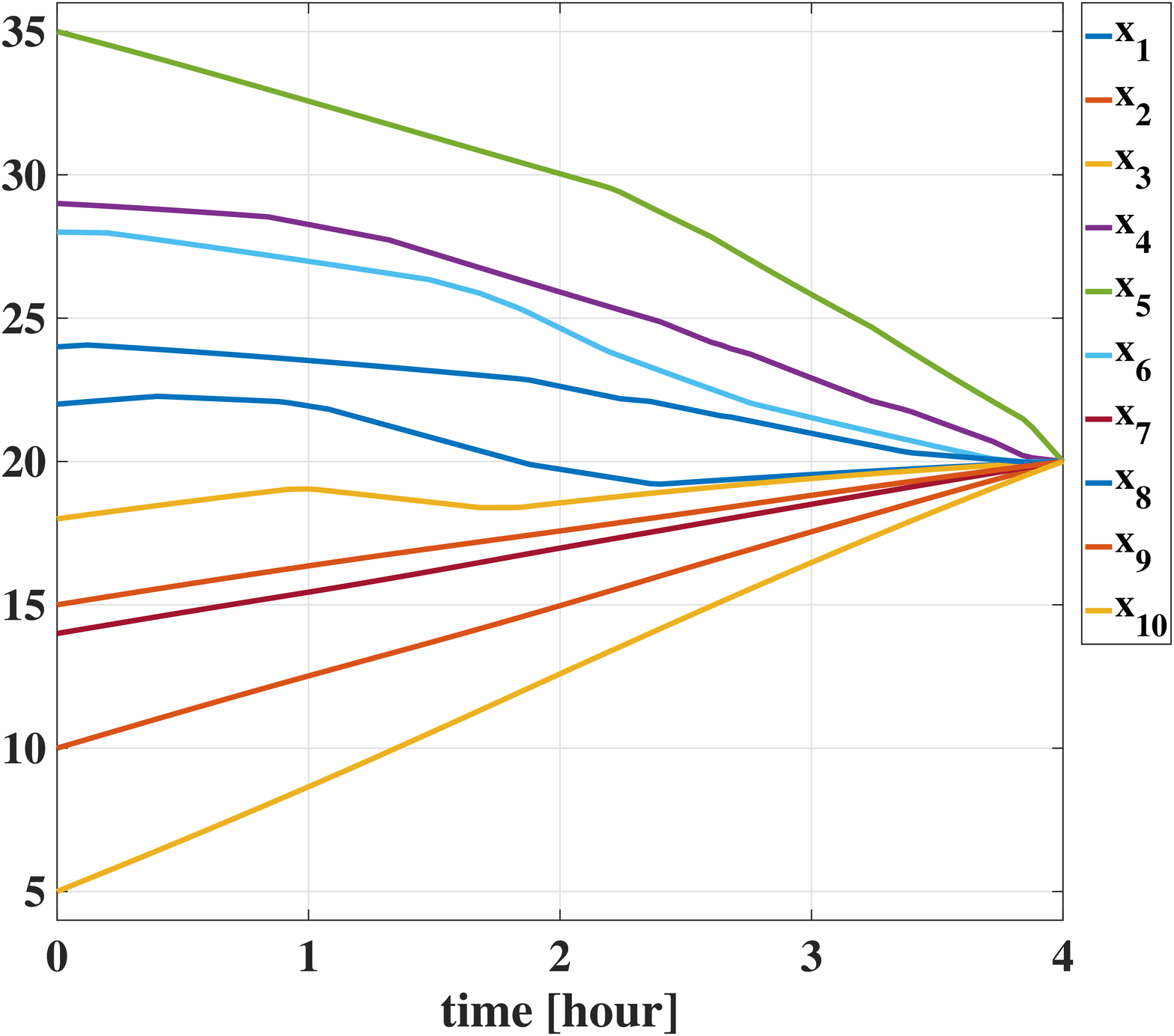}
    \includegraphics[width=0.75\linewidth]{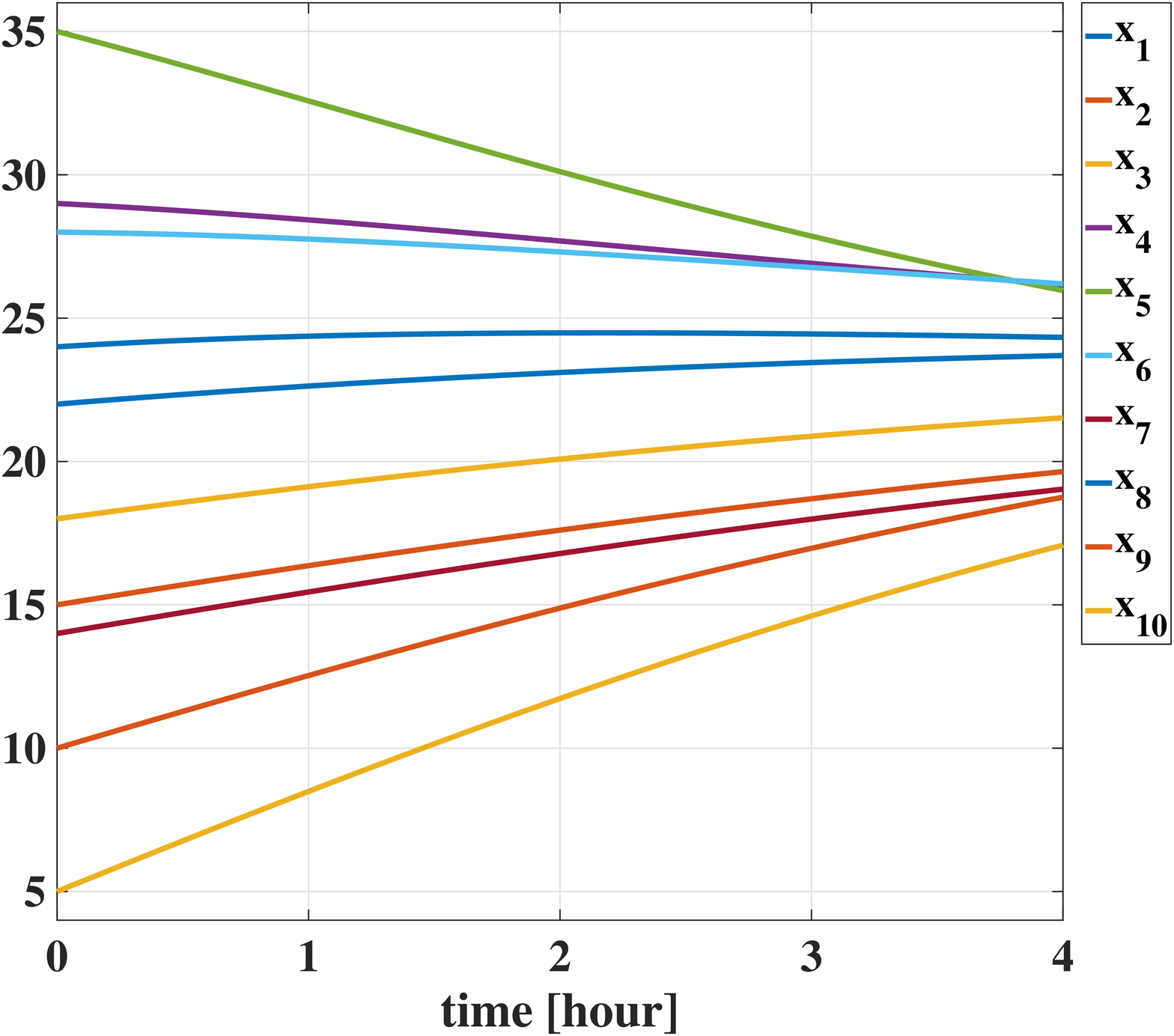}
  \caption{\red{State trajectories corresponding to the $L^0/\ell^0$ optimal control (top) and the zero control (bottom)}}
\label{fig:state_ex2}
\end{figure}


\begin{figure}[tbh]
  \centering
    \includegraphics[width=0.75\linewidth]{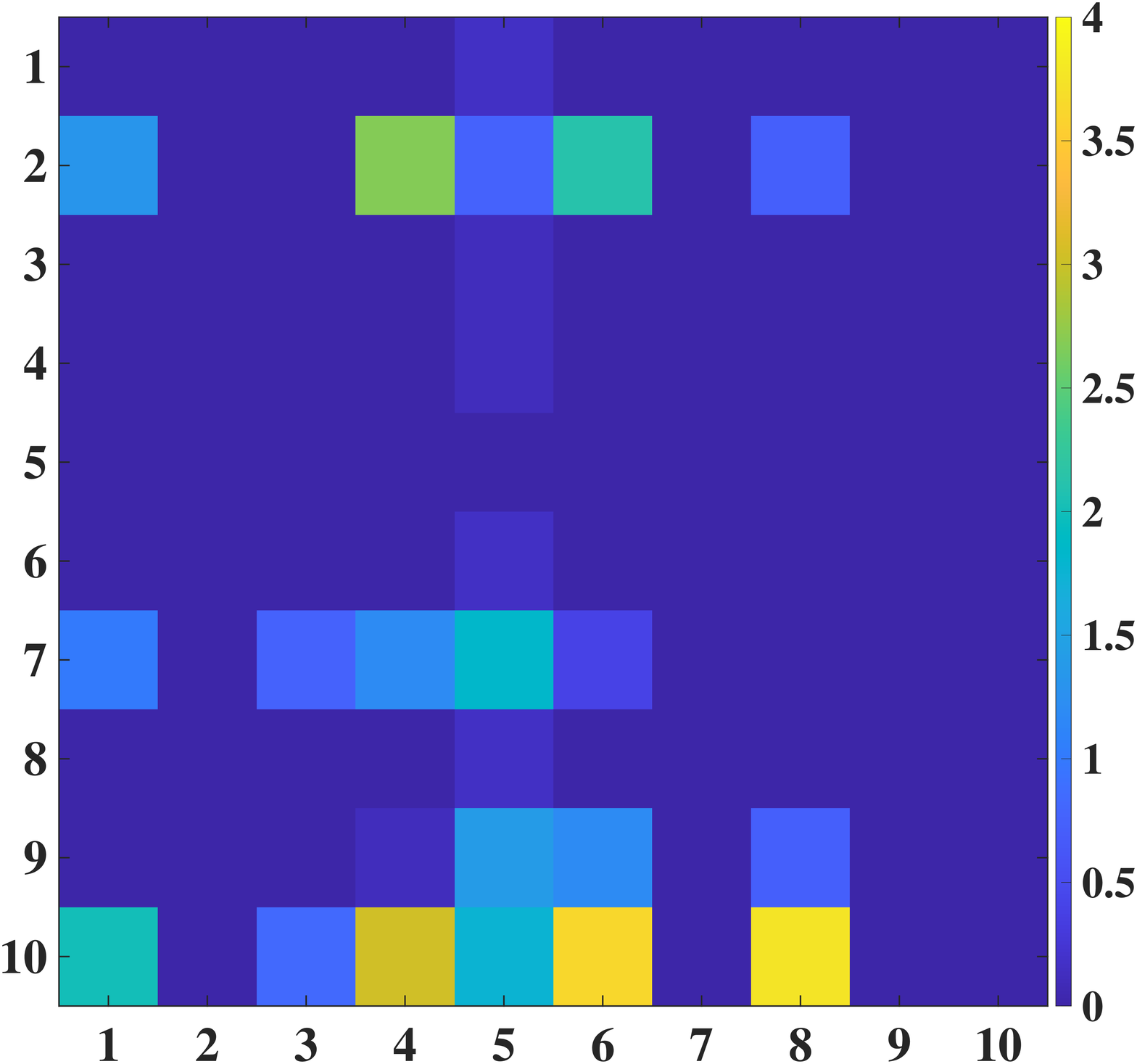}
    \includegraphics[width=0.75\linewidth]{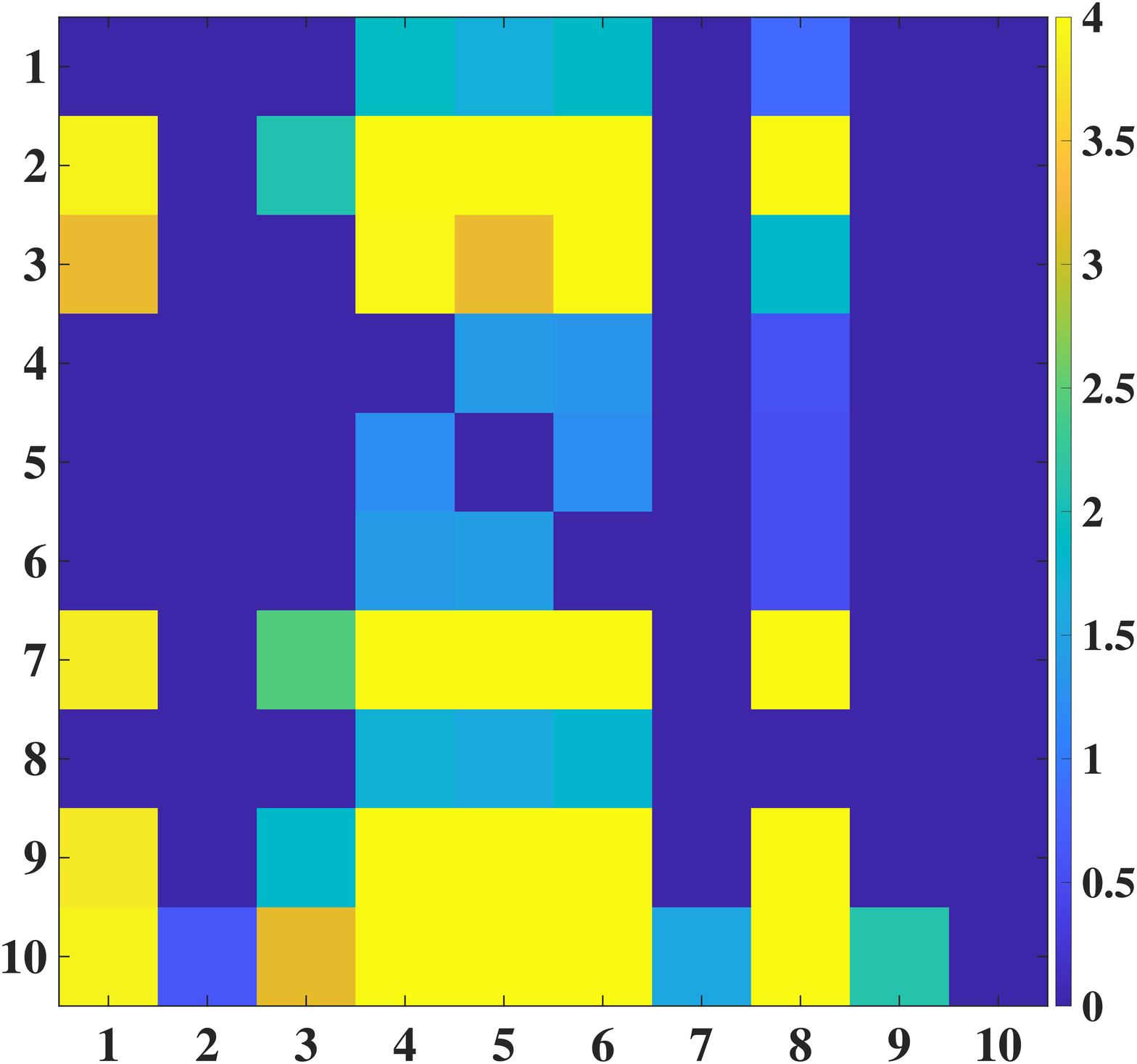}
  \caption{\red{Activated time duration of the trajectories obtained by the $L^0/\ell^0$ optimization (top) and the $L^2/\ell^1$ optimization (bottom)}}
\label{fig:control_cost_ex2a}
\end{figure}

We demonstrate the effectiveness of the developed method through numerical study. 
Let $s = 10$ be the number of the stations.
The total number of the vehicles is $200$,
which implies that $\sum_{i \in \mathcal{S}} x_i(t) = 200$ for all $t$.
The location (within a radius of about twenty kilometers) 
and initial number of vehicles of each station are depicted in Fig.~\ref{fig:station_ex2}. 
The color on the bottom indicates the degree of the road congestion in each place.
The desired final state $x_d$ is such that all the components are equally 20 for rebalancing the vehicles.
The price elasticity coefficients $\theta_{ij}(t)$ and the price adjustment coefficients $\lambda_{ij}(t)$ 
are taken from the uniform distribution on $[0, 0.3]$.
The ratio $\gamma_{ij}(t)$ is determined depending on the congestion and distance between the stations described in Fig. \ref{fig:station_ex2}.
This system is modeled as (\ref{eq_sk:dot_x})
with $A(t)$ and $B(t)$ given in (\ref{eq_sk:AB}).
The number of the service staff is set to $\beta=10$.
For this system, the proposed $L^0/\ell^0$ optimization method is applied for sparse rebalance.

Figs. \ref{fig:state_ex2} and \ref{fig:control_cost_ex2a} show the simulation results under this setting.
The upper panel of Fig. \ref{fig:state_ex2} represents the time plots of 
the components $x_i(t)$, the numbers of the vehicles in the stations, which 
shows that all $x_i(t)$ converge to 20 to achieve the rebalance after $4$ hours.
On the other hand, the lower panel of Fig. \ref{fig:state_ex2} represents those when no control input is applied, namely, $u_{ij}(t) = 0$,
where the components $x_i(t)$ do not agree at that time.
These results show that the proposed method accelerates the rebalance speed.
The upper panel of Fig. \ref{fig:control_cost_ex2a} represents the activated time duration ($L^0$ norm) of the control input $u_{ij}(t)$ when the proposed $L^0/\ell^0$ optimization method is applied,
while the lower panel shows that by using an $L^2/\ell^1$ optimization method.
These figures illustrate that many of the control inputs are equal to zero (inactive) via the $L^0/\ell^0$ optimization 
while many of the control inputs via the $L^2/\ell^1$ optimization are non-zero (active).
\red{Indeed, the $L^0$ costs of the optimal controls are approximately $32.08$ ($L^0/\ell^0$ optimization) and $136.36$ ($L^2/\ell^1$ optimization), respectively.}
Hence, sparse rebalance is successful due to the proposed method.

\section{Conclusion}
\label{sec:conclusion}
This paper has analyzed an optimal node scheduling that maximizes the trace of the controllability Gramian.
This analysis enables us to find an activation schedule of control inputs 
that steers the system while saving energy. 
Taking the number of control nodes and the time length of providing inputs into account,
our optimization problem newly includes two types of constraints on sparsity.
We have shown a sufficient condition under which
our sparse optimization problem boils down to a convex optimization problem.
This paper assumes the network topology among nodes is given and fixed.
Future work includes the design of the time-varying topology
and \red{more practical model of the staff dynamics in the mobility system}.


\bibliographystyle{IEEEtran}
\bibliography{main}
\end{document}


\maketitle

\section{Model derivation of the traffic system}

Assume that the customers take the service according to the Poisson process,
and that the time it takes to travel between stations follows an exponential distribution.
The latter assumption is valid because the farther customers are away from their destinations, the more they have chance of wasting time by shopping or getting lost on the way.
Then, the traffic system of customers' vehicles is modeled as
(17) and (18) in the main paper, that is,
\begin{align}
\dot f_{ij}(t) &=  -\gamma_{ij}(t) f_{ij}(t) + g_{ij}(t) + u_{ij}(t)
 \tag{r1}
 \\
 \dot v_i(t) &=  \sum_{j=1}^s \gamma_{ij}(t) f_{ij}(t)-\sum_{j=1}^s (g_{ji}(t) + u_{ji}(t))
 \tag{r2}
 .
\end{align}
See Table I for the meaning and units of the variables,
where [units] indicates the unit of the number of vehicles.

\begin{table}[bt]
 \centering
 \caption{Meaning and units of the variables}
  \begin{tabular}{ccc}
   \hline
   Variable & Unit & Meaning \\
   \hline
   $f_{ij}(t)$ & units & Expectation of the number of vehicles traveling from station $j$ to $i$\\
   $v_i(t)$ & units & Expectation of the number of vehicles parking at station $i$\\
   $g_{ij}(t)$ & units$\cdot$time$^{-1}$ & Number of the demand for travels from station $j$ to $i$ per time (the rate of the Poisson process)\\
   $u_{ij}(t)$ & units$\cdot$time$^{-1}$ &
    Number of rebalanced vehicles from station $j$ to $i$ per time\\
   $\gamma_{ij}(t)$ & time$^{-1}$ & Arrival rate of vehicles form station $j$ to $i$ per time (the inverse of $\tau_{ij}(t)$)\\
   $\tau_{ij}(t)$ & time & Average of the time it takes to travel from station $j$ to $i$\\
   \hline
  \end{tabular}
\end{table}

We derive the system consisting in (r1) and (r2) from the assumptions.
Here, the following notation is used.
For a random variable $X(t) \in \mathbb{Z}$, $\mathrm{P}[X(t) = k]$ denotes the probability that $X(t) = k$ occurs.
The expectation of $X(t)$ is represented by $\mathrm{E}[X(t)] := \sum_{k \in \mathbb{Z}} k\mathrm{P}[X(t) = k]$.

Let $F_{ij}(t) \in \mathbb{Z}$ [units] be the number of vehicles traveling from station $j$ to $i$ at time $t$,
and let $V_i(t) \in \mathbb{Z}$ [units] be the number of vehicles parking at station $i$ at time $t$.
These numbers vary within time interval $[t,t+\delta)$
for $\delta > 0$ according to
\begin{align}
 F_{ij}(t+\delta) &=  F_{ij}(t) + G_{ij}(t,\delta) + U_{ij}(t,\delta) - A_{ij}(t,\delta)
 \tag{r3}
 \\
 V_i(t+\delta) &=  V_i(t) +  \sum_{j=1}^s A_{ij}(t,\delta) - \sum_{j=1}^s (G_{ji}(t,\delta) + U_{ji}(t,\delta))
 \tag{r4}
,
\end{align}
where $G_{ij}(t,\delta) \in \mathbb{Z}$ [units] is the number of vehicles of customers departing from station $j$ for the destination $i$,
$U_{ij}(t,\delta) \in \mathbb{Z}$ [units] is the number of vehicles for rebalancing, departing from station $j$ for the destination $i$,
and $A_{ij}(t,\delta) \in \mathbb{Z}$ [units] is the number of vehicles arriving at station $i$ from $j$ within interval $[t,t+\delta)$.

We assume the following conditions on
$G_{ij}(t,\delta)$, $U_{ij}(t,\delta)$, and
$A_{ij}(t,\delta)$.
First, the number $G_{ij}(t,\delta)$ of the customers is a random variable following the Poisson process
\begin{equation}
 \mathrm{P}[G_{ij}(t,\delta) = k] = \frac{1}{k!}(g_{ij}(t)\delta)^k e^{-g_{ij}(t)\delta}
 \tag{r5}
\end{equation}
for $k \in \mathbb{Z}$ with a rate $g_{ij}(t) \ge 0$ [units$\cdot$time$^{-1}$],
which is the number of the demand for traveling from station $j$ to $i$ per time.
From (r5),
the expectation of $G_{ij}(t,\delta)$ is calculated as follows:
\begin{equation}
  \mathrm{E}[G_{ij}(t,\delta)] = g_{ij}(t)\delta
  .
  \tag{r6}
\end{equation}
Second, the number $U_{ij}(t,\delta)$ of rebalanced vehicles
from station $j$ to $i$ is determined
with a rate $u_{ij}(t) \ge 0$ [units$\cdot$time$^{-1}$],
which leads to
\begin{equation}
 \mathrm{E}[U_{ij}(t,\delta)] = u_{ij}(t)\delta
 .
 \tag{r7}
\end{equation}
Third, let $p_{ij}(t,\delta) \in (0,1)$ be the the probability of the arrival
of the vehicle traveling from station $j$ to $i$
within interval $[t,t+\delta)$.
Then, the number $A_{ij}(t,\delta)$ of vehicles arriving at station $i$ from $j$
within interval $[t,t+\delta)$ is given
as the random variable following the binomial distribution
\begin{equation}
  \mathrm{P}[A_{ij}(t,\delta) = k|F_{ij}(t)=\ell] = {}_{\ell}\mathrm{C}_k p_{ij}(t,\delta)^k (1-p_{ij}(t,\delta))^{\ell-k}
 \tag{r8}
\end{equation}
for $k,\ell \in \mathbb{Z}$,
which indicates the probability that $k$ vehicles arrive at station $i$ out of $\ell$ vehicles.
From (r8),
the conditional expectation of $A_{ij}(t,\delta)$ is calculated as
\begin{equation}
\mathrm{E}[A_{ij}(t,\delta)|F_{ij}(t)] = p_{ij}(t,\delta) F_{ij}(t)
 \tag{r9}
.
\end{equation}
We assume that the time $T_{ij}(t) > 0$ it takes to travel from $j$ to $i$ is the random variable following an exponential distribution
with average $\tau_{ij}(t) > 0$.
Then, the probability $p_{ij}(t,\delta)$ of the arrival is of the form
\begin{equation}
 p_{ij}(t,\delta) =
 \mathrm{P}[T_{ij}(t) < \delta] = 1 - e^{-\frac{\delta}{\tau_{ij}(t)}}
 \tag{r10}
 .
\end{equation}
%
From (r9), (r10), and the property of the conditional expectation,
the expectation of $A_{ij}(t,\delta)$ is calculated as
\begin{align}
\mathrm{E}[A_{ij}(t,\delta)] &=
\mathrm{E}[\mathrm{E}[A_{ij}(t,\delta)|F_{ij}(t)]]
\notag\\
&=\mathrm{E}[p_{ij}(t,\delta)F_{ij}(t)]
= (1-e^{-\frac{\delta}{\tau_{ij}(t)}})\mathrm{E}[F_{ij}(t)]
\tag{r11}
.
\end{align}

Let $f_{ij}(t) = \mathrm{E}[F_{ij}(t)]$ [units]
be the expectation of the number of vehicles
traveling from $j$ to $i$,
and let $v_i(t) = \mathrm{E}[V_i(t)]$ [units]
be the expectation of the number of vehicles
parking at $i$.
Take the expectations of (r3) and (r4),
and from (r6), (r7), and (r11),
we obtain
\begin{align}
 f_{ij}(t+\delta) -f_{ij}(t) &=   g_{ij}(t)\delta + u_{ij}(t)\delta - (1-e^{-\frac{\delta}{\tau_{ij}(t)}})f_{ij}(t)
 \tag{r12}
 \\
 v_i(t+\delta) - v_i(t) &= \sum_{j=1}^s (1-e^{-\frac{\delta}{\tau_{ij}(t)}})f_{ij}(t)
 - \sum_{j=1}^s (g_{ji}(t)\delta + u_{ji}(t)\delta)
 \tag{r13}
.
\end{align}
Divide (r12) and (r13) by $\delta$ and take the limit
$\delta \to 0$.
Then, we obtain (r1) and (r2)
for $\gamma_{ij}(t) = 1/\tau_{ij}(t)$ [time$^{-1}$],
which indicates the arrival rate of vehicles traveling from $j$ to $i$
per time.

\section{Model of rebalancing by staff}

In the main paper, rebalancing by staff is modeled only with static constraints.
This model is valid under assumption that the dynamics of the staff is sufficiently
faster than that of customers,
which is fairly realistic because the dynamics of customers includes not only traveling between stations by vehicles, but also reserving vehicles, walking to stations, and so forth.
However, if a more precise model is required, a dynamic model can be obtained
in the same way as (r1) and (r2) as follows:
\begin{align}
 \dot {\hat f}_{ij}(t) &=  {\hat u}_{ij}(t) - \hat\gamma_{ij}(t) \hat f_{ij}(t)
 \tag{r14}
 \\
 \dot{\hat v}_i(t) &=  \sum_{j=1}^s \hat\gamma_{ij}(t) \hat f_{ij}(t)-\sum_{j=1}^s \hat u_{ji}(t)
 \tag{r15}
,
\end{align}
where the variables with $\hat{~}$ refer to the variables concerning the staff
in the same meaning as Table I.
Assume that each team of the staff drives their own car
with/without accompanying a vehicle of customers for rebalance,
which leads to the constraint
\begin{equation}
 \hat u_{ij}(t) \ge u_{ij}(t)
 \tag{r16}
 .
\end{equation}
In (r16), if a vehicle of customers is accompanied, the equality holds;
otherwise, it does not.
Then, the control inputs of the whole system
are $u_{ij}(t)$ and $\hat u_{ij}(t)$ satisfying (r16).
Furthermore, the non-negativeness of
the number of the staff cars is required as
\begin{equation}
  \hat v_i(t) \ge 0
  \tag{r17}
  .
\end{equation}

Next, we consider the constraints on the individual staff teams.
There are $\beta$ teams of the staff for rebalance,
and let $u_{ij,k}(t) \ge 0$ be the number of the vehicles
rebalanced by team $k \in \{1,2,\ldots,\beta\}$.
Each team can move one vehicle on one route at once,
which is expressed as
\begin{equation}
 \|u_k\|_{L^\infty} \le 1,~\|u_k(t)\|_{\ell^0} \le 1
 \tag{r18}
 ,
\end{equation}
where $u_k = [u_{12,k}~ u_{13,k}~\cdots u_{n(n-1),k}]^\top$.
Note that the $\ell^0$ norm plays an essential role in (r18)
because one team cannot move more than one vehicle simultaneously.
The total number of the vehicles for rebalance by all the teams
from station $j$ to $i$ is given as
\begin{equation}
 u_{ij}(t) = \sum_{k=1}^\beta u_{ij,k}(t)
 \tag{r19}
.
\end{equation}
Replace $u_{ij}(t)$ with (r19) in (r1) and (r2),
and we obtain the traffic model with the control input $u_{ij,k}(t)$.

Combining the above two models yields the model of the staff considering the dynamics
and constraints of individual teams as follows:
\begin{align}
 \dot {\hat f}_{ij,k}(t) &=  {\hat u}_{ij,k}(t) - \hat\gamma_{ij,k}(t) \hat f_{ij,k}(t)
 \tag{r14'}
 ,\\
 \dot{\hat v}_{i,k}(t) &=  \sum_{j=1}^s \hat\gamma_{ij,k}(t) \hat f_{ij,k}(t)-\sum_{j=1}^s \hat u_{ji,k}(t)
 \tag{r15'}
,\\
 \hat u_{ij,k}(t) &\ge u_{ij,k}(t)
 \tag{r16'}
,\\
 \hat v_{i,k}(t) & \ge 0
 \tag{r17'}
 ,\\
 \|\hat u_k\|_{L^\infty} &\le 1,~\|\hat u_k(t)\|_{\ell^0} \le 1
 \tag{r18'}
 ,\\
  u_{ij}(t) &= \sum_{k=1}^\beta u_{ij,k}(t)
 \tag{r19'}
 .
\end{align}
The optimization problem with this model is formulated as
\begin{equation*}
 \underset{u,\hat u}{\text{minimize}}~~ \|u\|_{L^0} + \|\hat u\|_{L^0}
\end{equation*}
under the constraints (r1), (r2), and (r14')--(r19')
for the variables
$u = [u_1~ \cdots~u_\beta]^\top$ and
$\hat u = [\hat u_1~ ~ \cdots~\hat u_\beta]^\top$,
where $u_k = [u_{12,k}~u_{13,k}~\cdots~u_{n(n-1),k}]^\top$
and $\hat u_k = [\hat u_{12,k}~\hat u_{13,k}~\cdots~\hat u_{n(n-1),k}]^\top$.
\red{This problem can be written in the same form as Problem 3 in the main part. Therefore, the relaxed convex problem, for which $L^0$ and $\ell^0$-norms are replaced by $L^1$ and  $\ell^1$-norms, can be numerically solved via time-discretization. Note, however, that (r16') is not a box constraint, and that (r17') is a state constraint.
The theoretical equivalence between the original and relaxed problems under these constraints are beyond the focus of this paper, and left as future work. }
